\documentclass[journal, onecolumn]{IEEEtran}
\usepackage{cite}

\ifCLASSINFOpdf
   \usepackage[pdftex]{graphicx}
\else
   \usepackage[dvipdfmx]{graphicx}
\fi
\usepackage{amsmath,amsmath,amsfonts,amssymb}
\newtheorem{thm}{Theorem}[section]
\newtheorem{corollary}{Corollary}[section]
\newtheorem{lemma}{Lemma}[section]
\newtheorem{prop}{Proposition}[section]

\newcommand{\tr}{\mathrm{tr}}

\newcommand{\KL}[2]{D_{\mathrm{KL}} \left( #1 \,\middle|\!\middle|\, #2 \right)}

\newcommand{\pd}[1]{\partial_{#1}}
\newcommand{\pu}[1]{\partial^{#1}}
\newcommand{\pHd}[1]{\partial_{#1}}
\newcommand{\pAd}[1]{\partial_{\bar{#1}}}
\newcommand{\pHu}[1]{\partial^{#1}}
\newcommand{\pAu}[1]{\partial^{\bar{#1}}}

\newcommand{\gdd}[2]{g_{#1 #2}}
\newcommand{\guu}[2]{g^{#1 #2}}

\newcommand{\gHdHd}[2]{g_{#1 #2}}
\newcommand{\gAdAd}[2]{g_{\bar{#1} \bar{#2}}}
\newcommand{\gHdAd}[2]{g_{#1 \bar{#2}}}
\newcommand{\gAdHd}[2]{g_{\bar{#1} #2}}

\newcommand{\gHuHu}[2]{g^{#1 #2}}
\newcommand{\gAuAu}[2]{g^{\bar{#1} \bar{#2}}}
\newcommand{\gHuAu}[2]{g^{#1 \bar{#2}}}
\newcommand{\gAuHu}[2]{g^{\bar{#1} #2}}

\newcommand{\Guddd}[4]{{\overset{\,\tiny{#1}}{\Gamma}}_{#2 #3 #4}}
\newcommand{\Guddu}[4]{ \overset{\tiny{#1} \quad\quad }{ {{\Gamma}_{#2 #3}}^{#4} } }

\newcommand{\GuHd}[2]{ \overset{\tiny{#1} \quad }{ {\Gamma}_{#2} }}

\newcommand{\GuHdHdAd}[4]{ \overset{\tiny{#1} \quad }{ {{\Gamma}_{#2 #3 \bar{#4}}} } }

\newcommand{\GuHdHdHu}[4]{ \overset{\tiny{#1} \quad }{ {{\Gamma}_{#2 #3}}^{#4} } }
\newcommand{\GuHdAdAu}[4]{ \overset{\tiny{#1} \quad }{ {{\Gamma}_{#2 \bar{#3}}}^{\bar{#4}} } }
\newcommand{\GuAdAdAu}[4]{ \overset{\tiny{#1} \quad }{ {{\Gamma}_{\bar{#2} \bar{#3}}}^{\bar{#4}} } }

\newcommand{\Tddd}[3]{T_{#1 #2 #3}}
\newcommand{\THdHdAd}[3]{T_{#1 #2 \bar{#3}}}
\newcommand{\THdAdHd}[3]{T_{#1 \bar{#2} #3}}

\newcommand{\Td}[1]{T_{#1}}
\newcommand{\THd}[1]{T_{#1}}

\begin{document}

\title{Shrinkage priors on complex-valued circular-symmetric autoregressive processes}
\author{Hidemasa~Oda
and~Fumiyasu~Komaki%
\thanks{H. Oda is with the Department of Mathematical Informatics,
The University of Tokyo, Tokyo 113-8654, Japan.}%
\thanks{F. Komaki is with the Department of Mathematical Informatics,
The University of Tokyo, Tokyo 113-8654, Japan.}
}

\markboth{}%
{Oda and Komaki: Shrinkage priors on complex-valued circular-symmetric autoregressive processes}

\maketitle

\begin{abstract}
We investigate shrinkage priors on power spectral densities for complex-valued circular-symmetric autoregressive processes.
We construct shrinkage predictive power spectral densities, which asymptotically dominate (i) the Bayesian predictive power spectral density based on the Jeffreys prior and (ii) the estimative power spectral density with the maximum likelihood estimator,
where the Kullback--Leibler divergence from the true power spectral density to a predictive power spectral density is adopted as a risk.
Furthermore, we propose general constructions of objective priors for K\"ahler parameter spaces
by utilizing a positive continuous eigenfunction of the Laplace--Beltrami operator with a negative eigenvalue.
We present numerical experiments on a complex-valued stationary autoregressive model of order 1. 
\end{abstract}


\IEEEpeerreviewmaketitle

\section{Introduction}

We investigate the time fluctuation of a single particle having a circular-symmetric distribution in a two-dimensional space.
In this situation, the complex plane $\mathbb{C}$ is often used for the representation of the process,
and an observation of a single particle at different time points can be represented as a complex-valued vector.
A complex-valued random vector $Z$ is called circular-symmetric if, for any constant $\phi \in \mathbb{R}$, the distribution of $e^{\sqrt{-1} \, \phi} Z$ equals the distribution of $Z$.

We focus on complex-valued circular-symmetric discrete Gaussian processes, which are defined as complex-valued processes having finite-dimensional marginal distributions that are complex normal distributions.
The precise definitions of a complex normal distribution and a complex-valued Gaussian process are given in Section \ref{section:complex_valued_gaussian_process}.
The circular symmetry of complex normal distributions with zero mean is of great importance in practical applications \cite{miller1974complex}.
Complex-valued processes are commonly used for directional processes, such as wind, radar, and sonar signals \cite{mandic2009complex}.
Furthermore, complex-valued representations are widely used in diverse fields, such as econometrics \cite{svetunkov2012complex} and complex-valued neural networks \cite{hirose2003complex}.

The power spectral density of a complex-valued Gaussian process is defined as a Fourier transform
\begin{align}
  S(\omega) := \frac{1}{2\pi} \sum_{h = -\infty}^{\infty} \gamma_h \, e^{-\sqrt{-1} h \omega}
  \label{definition:power_spectral_density}
\end{align}
of the autocovariances $\{ \gamma_h \}_{h \in \mathbb{Z}}$ of the process, where $\omega \in [-\pi, \pi]$.
We parametrize complex-valued Gaussian processes by complex variables $\theta = (\theta^1, \cdots, \theta^p)$ $\in \Theta \subset \mathbb{C}^p$.
In other words, we regard the sequence $\{ \gamma_h \}_{h \in \mathbb{Z}}$ of autocovariances of the process as functions $\{ \gamma_h (\theta) \}_{h \in \mathbb{Z}}$ of the parameter $\theta$.
For each $\theta \in \Theta$, we denote the corresponding power spectral density by $S_{\theta} (\omega)$ or $S (\omega \mid \theta)$.

We aim to predict the joint distribution of future observations $w = \left( w_1, w_2, \cdots \right)$ by observing current observations $z^{(N)} = (z_1, \cdots, z_N) \in \mathbb{C}^N$ of size $N$ from a complex-valued Gaussian process having the true parameter $\theta_0 \in \Theta$.
We suppose that the sample $w$ and the sample $z^{(N)}$ are independent;
i.e., the sample $w$ is taken from a different process of the same type or from the same process but a long time after the sample $z^{(N)}$ is taken.
Let us consider the problem of constructing a power spectral density $\hat{S}^{(N)}$ that corresponds to the joint distribution of future observations $w$; see \cite{komaki1999estimating}.
The constructed power spectral density $\hat{S}^{(N)}$ is called a predictive power spectral density.
More precisely, a predictive power spectral density $\hat{S}^{(N)} (\omega)$ is a function of an observation $z^{(N)}$ for each $\omega \in [-\pi, \pi]$.
The goodness of the prediction is evaluated by the risk, which is defined as
\begin{align}
  R \bigl( \hat{S}^{(N)} \mid \theta_0 \bigr)
  &:= E_{\theta_0} \left[ \KL{S_{\theta_0}}{\hat{S}^{(N)}} \right] 
  = \int_{\mathbb{C}^N} \KL{S_{\theta_0}}{\hat{S}^{(N)}} \,dP^{(N)}_{\theta_0} \bigl( z^{(N)} \bigr)
  ,\label{definition:risk}
\end{align}
where
$P^{(N)}_{\theta_0}$ denotes the distribution of $z^{(N)}$ and
the Kullback--Leibler divergence between two power spectral densities $S_1$ and $S_2$ is defined as
\begin{align}
  \KL{S_1}{S_2} &:= \frac{1}{2\pi} \int_{-\pi}^{\pi} \left( - \log \frac{S_1 (\omega)}{S_2 (\omega)} - 1 + \frac{S_1 (\omega)}{S_2 (\omega)} \right) d\omega
  \label{definition:kl_divergence}
  \,;
\end{align}
see Appendix \ref{section:kl_divergence}.
The principal aim of this study is to construct a predictive power spectral density $\hat{S}^{(N)}$ with its risk $R \bigl( \hat{S}^{(N)} \mid \theta \bigr)$ being as small as possible for most $\theta \in \Theta$.
There is another interpretation for the problem:
we estimate the underlying true power spectral density of the process under the loss given by (\ref{definition:risk}) based on the Kullback--Leibler divergence.

\IEEEpubidadjcol

There are two basic constructions for predictive power spectral densities.
The first construction, called the estimative method, is $S_{\hat{\theta}^{(N)}}$, where $\hat{\theta}^{(N)} = \hat{\theta}^{(N)} \bigl( z^{(N)} \bigr)$ is an estimator for the true parameter $\theta_0$.
The second construction, called the Bayesian predictive method, is
\begin{align}
  \hat{S}^{(N)}_{\pi} (\omega) := \int_{\Theta} S \bigl( \omega \mid \theta \bigr) \, \pi \bigl( \theta \mid z^{(N)} \bigr) \,d\theta
  \label{definition:BPD}
\end{align}
for a possibly improper prior $\pi$.
Here, 
\begin{align}
  \pi \bigl( \theta \mid z^{(N)} \bigr) := \frac{ p_{\theta} \bigl( z^{(N)} \bigr) \pi (\theta) }{ \int_{\Theta} p_{\theta'} \bigl( z^{(N)} \bigr) \pi (\theta') \,d\theta' }
\end{align}
denotes the posterior, given an observation $z^{(N)} \in \mathbb{C}^N$ based on the prior $\pi$, and $\int_{\Theta} \pi \bigl( \theta \mid z^{(N)} \bigr) \,d\theta < \infty$ is assumed.
The power spectral density $\hat{S}^{(N)}_{\pi}$ is called the Bayesian predictive power spectral density based on the prior $\pi$.
Herglotz's theorem asserts that a family $\{ \gamma_h \}_{h \in \mathbb{Z}}$ of complex numbers parameterized by integers is the family of autocovariances of a complex-valued stationary process if and only if there exists a spectral distribution function $F$ such that $\gamma_h = \int_{-\pi}^{\pi} e^{\sqrt{-1} h \omega} \,dF(\omega)$ and $\{ \gamma_h \}_{h \in \mathbb{Z}}$ is positive semi-definite; i.e., $\sum_{i=1}^n \sum_{j=1}^n a_i \, \gamma_{i-j} \, \overline{a}_j \geq 0$ for any $n > 0$ and $(a_1, \cdots, a_n) \in \mathbb{C}^n$; see Corollary 4.3.1 in \cite{brockwell1991time}.
Because
\begin{align}
  \hat{\gamma}^{(N)}_{\pi, h}
  := \int_{-\pi}^{\pi} e^{\sqrt{-1} h \omega} \, \hat{S}^{(N)}_{\pi} ( \omega ) \,d\omega
  = \int_{\Theta} \gamma_h (\theta) \, \pi \bigl( \theta \mid z^{(N)} \bigr) \,d\theta
\end{align}
is a weighted average of $\gamma_h (\theta)$,
the process with the predictive power spectral density given by (\ref{definition:BPD}) is stationary as long as each process parameterized by $\theta \in \Theta$ is stationary.
If a prior $\pi$ is given and the risk is defined as (\ref{definition:risk}), the Bayesian predictive power spectral density $\hat{S}^{(N)}_{\pi}$ minimizes the Bayes risk
\begin{align} 
  r \bigl( \hat{S}^{(N)} \mid \pi \bigr) := \int_{\Theta} R \bigl( \hat{S}^{(N)} \mid \theta \bigr) \pi(\theta) \,d\theta
  \label{definition:bayes_risk}
\end{align}
among all the predictive power spectral densities $\hat{S}^{(N)}$ as long as the Bayes risk $r \bigl( \hat{S}_{\pi}^{(N)} \mid \pi \bigr)$ is finite.
Therefore, the remaining problem is to determine and construct an appropriate prior $\pi$.

Non-informative priors for time-series models, such as the Jeffreys prior, which is usually improper, have been discussed in previous works \cite{berger1994noninformative, zellner1977maximal, philippe2002non}.
We propose a proper prior $\pi^{(-1)}$ defined on the complex parameter space $\Theta \subset \mathbb{C}^p$ for the complex-valued stationary autoregressive processes $\mathrm{AR} (p; \mathbb{C})$ of order $p \geq 1$.
The Bayesian predictive power spectral density $\hat{S}^{N}_{\pi^{(-1)}}$ based on the proposed prior $\pi^{(-1)}$ asymptotically dominates the estimative power spectral density $S_{\hat{\theta}^{(N)}}$ with the maximum likelihood estimator $\hat{\theta}^{(N)}$.
Moreover, the proposed predictive power spectral density $\hat{S}^{(N)}_{\pi^{(-1)}}$ asymptotically dominates the Bayesian predictive power spectral density $\hat{S}^{(N)}_{\pi_J}$ based on the Jeffreys prior $\pi_J$,
and the $O(N^{-2})$ term of the risk improvement is constant regardless of $\theta \in \Theta$:
\begin{align}
  R \bigl( \hat{S}^{(N)}_{\pi_J} \mid \theta \bigr) - R \bigl( \hat{S}^{(N)}_{\pi^{(-1)}} \mid \theta \bigr)
  = \frac{2 p (p + 1)}{N^2}  + O(N^{-\frac{5}{2}})
  , \label{formula:main_ar}
\end{align}
which is summarized as the Main Theorem in Section \ref{section:main_theorem}.

An eigenfunction $\phi$ of the Laplacian (Laplace--Beltrami operator) $\Delta$ plays a crucial role in constructing the proposed prior $\pi^{(-1)}$.
The Laplacian is a differential operator that does not depend on the specific choice of parameterizations of the parameter space; see Appendix \ref{section:wirtinger_calculus}.
This operator transforms a scalar function defined on the parameter space to another scalar function defined on the parameter space.
References \cite{brown1971admissible} and \cite{stein1981estimation} stress the importance of the super-harmonicity for the shrinkage effect in the estimation of the mean of a multivariate normal distribution.
More generally, it is known that the super-harmonicity of the ratio of the proposed prior to the Jeffreys prior is the key to inducing the shrinkage effect \cite{komaki2006shrinkage}.
Another important property in the construction of the proposed prior $\pi^{(-1)}$ is the K\"ahlerness, the generalization of the concept of exponential families, of the complex parameter space $\Theta$.
The parameter space of the complex-valued stationary autoregressive moving average processes $\mathrm{ARMA} (p, q; \mathbb{C})$ is shown to be K\"ahler in \cite{choi2015kahlerian}.
We give a general construction of priors utilizing a positive continuous eigenfunction $\phi > 0$ of the Laplacian $\Delta$ with a negative eigenvalue $-K < 0$, i.e., $\Delta \phi = -K \phi < 0$.
We define a family of priors $\{ \pi^{(\kappa)} \}_{\kappa \in \mathbb{R}}$ by $\pi^{(\kappa)} := \phi^{-\kappa+1} \pi_J$, which are called $\kappa$-priors.
We prove that if $-1 \leq \kappa < 1$, then $\hat{S}^{(N)}_{\pi^{(\kappa)}}$ asymptotically dominates $\hat{S}^{(N)}_{\pi_J}$.
To maximize the worst case of the risk improvement, we propose the $\kappa$-prior for $\kappa = -1$ that achieves a constant risk improvement.

The remainder of this paper is organized as follows. 
In Section \ref{section:complex_valued_gaussian_process}, we give the asymptotic expansion of Bayesian predictive power spectral densities $\hat{S}^{(N)}_{\pi}$ for complex-valued autoregressive moving average processes.
In Section \ref{section:complex_valued_arma_process}, a specific K\"ahler parameterization for $\mathrm{AR} (p; \mathbb{C})$ is introduced.
In Section \ref{section:main_theorem}, the Main Theorem (\ref{formula:main_ar}) is stated.
In Section \ref{section:super_harmonicity_prior}, we explicitly give the construction of the positive continuous eigenfunction $\phi$ with a negative eigenvalue $-K = -p(p+1)$ on the K\"ahler parameter space for $\mathrm{AR} (p; \mathbb{C})$.
Furthermore, we show that the family $\{ \pi^{(\kappa)} \}_{\kappa \in \mathbb{R}}$ of the proposed $\kappa$-priors $\pi^{(\kappa)}$ for $\mathrm{AR} (p; \mathbb{C})$ is a family of $\alpha$-parallel priors, which was introduced in \cite{takeuchi2005spl}, with $\alpha = (1 - \kappa) / 2$.
The generalization of the family of the proposed priors $\{ \pi^{(\kappa)} \}_{\kappa \in \mathbb{R}}$ and its relation with the $\alpha$-parallel priors for the i.i.d.\ case is discussed in Section \ref{section:generalization}.
In Section \ref{section:numerical_experiments}, numerical experiments are reported for the value of the risk differences $N^2 \bigl( R \bigl( \hat{S}^{(N)}_{\pi_J} \mid \theta \bigr) - R \bigl( \hat{S}^{(N)}_{\pi^{(-1)}} \mid \theta \bigr) \bigr)$ for $\mathrm{AR} (1; \mathbb{C})$.

\section{Bayesian predictive power spectral densities for complex-valued Gaussian processes}
\label{section:complex_valued_gaussian_process}

As explained in the Introduction, our aim is to construct a predictive power spectral density $\hat{S}^{(N)}$ after observing a sample $z^{(N)} \in \mathbb{C}^N$ of size $N$.
In this section, we present the asymptotic expansion of Bayesian predictive power spectral densities $\hat{S}^{(N)}_{\pi}$ for complex-valued autoregressive moving average processes.
This asymptotic expansion is a basic tool for assessing the performance of the choice of a prior $\pi$.
In the remainder of this paper, the Einstein notation is assumed.
Therefore, the summation is automatically taken over those indices that appear exactly twice, once as a superscript and once as a subscript.
The symbols $\alpha, \beta, \gamma, \cdots$ run through the indices $\{1, \cdots, p, \bar{1}, \cdots, \bar{p} \}$, and the symbols $i, j, k, \cdots$ run through the indices $\{ 1, \cdots, p \}$.

We first define the multivariate complex normal distribution.
Let $\mu \in \mathbb{C}^N$ and $\Sigma^{(N)}$ be an $N \times N$ complex-valued positive definite Hermitian matrix.
Note that the determinant $| \Sigma^{(N)} |$ of the matrix $\Sigma^{(N)}$ is positive.
An $N$-dimensional complex normal distribution (complex-valued circular-symmetric multivariate normal distribution) with mean $\mu$ and variance $\Sigma^{(N)}$
is defined by its probability density function:
\begin{align}
    p \left( z^{(N)} \mid \mu, \Sigma^{(N)} \right) := \frac{1}{\pi^N \, \left| \Sigma^{(N)} \right| } e^{- \, \left( z^{(N)} - \mu \right)^* \, \left( \Sigma^{(N)} \right)^{-1} \, \left( z^{(N)} - \mu \right)},
    \label{definition:complex_normal_distribution_probability_density_function}
\end{align}
where $z^{(N)} = (z^1, \cdots, z^N) \in \mathbb{C}^N$ and $z^{(N)*}$ denotes the complex conjugate transpose of $z^{(N)}$; see \cite{miller1974complex}.
The circular symmetry of a complex normal distribution with mean $0 \in \mathbb{C}^N$ is understood from the definition (\ref{definition:complex_normal_distribution_probability_density_function}).
The complex normal distribution with mean $0 \in \mathbb{C}^N$ and the identity matrix of size $N$ as the variance-covariance matrix is called the standard complex normal distribution of size $N$.
If we let $z^i = x^i + \sqrt{-1} \, y^i$ for $i = 1, \cdots, N$,
the $2N$-dimensional real-valued vector $(x^1, \cdots, x^N, y^1, \cdots, y^N)$ follows the $2N$-dimensional real-valued multivariate normal distribution with mean $(\Re (\mu), \Im (\mu))$ and the following variance-covariance matrix:
\begin{align}
  \begin{bmatrix}
      \frac{1}{2}\, \Re\, (\Sigma) & -\frac{1}{2}\, \Im\, (\Sigma) \\
      \frac{1}{2}\, \Im\, (\Sigma) & \frac{1}{2}\, \Re\, (\Sigma)
  \end{bmatrix}
  .
  \label{formula:viariance_covariance_matrix}
\end{align}
Therefore, an $N$-dimensional complex normal distribution is a special case of a $2N$-dimensional real normal distribution;
however, the opposite is not true.

A complex-valued discrete process $\{ Z_t \}_{t \in \mathbb{Z}}$ is called a Gaussian process (complex-valued circular-symmetric discrete Gaussian process) if
the tuple $(Z_{t_1}, Z_{t_2}, \cdots, Z_{t_N})$ of size $N$ follows a complex normal distribution for any $N$ and any $t_1, t_2, \cdots, t_N \in \mathbb{Z}$.
A complex Gaussian white noise $\{ \varepsilon_t \}_{t \in \mathbb{Z}}$ with variance $\sigma^2$ is a Gaussian process and $\frac{1}{\sigma} (\varepsilon_{t_1}, \cdots, \varepsilon_{t_N})$ follows a standard complex normal distribution of size $N$ for any $N$ and any $t_1, t_2, \cdots, t_N \in \mathbb{Z}$.
For a strongly stationary Gaussian process $\{ Z_t \}_{t \in \mathbb{Z}}$, we define the autocovariance $\gamma_h$ of order $h$ as the covariance of $Z_{t+h}$ and $\overline{Z_t}$.
Note that the autocovariances $\{ \gamma_h \}_{h \in \mathbb{Z}}$ are complex-valued, and the relation $\overline{\gamma_h} = \gamma_{-h}$ holds for any $h$.
The power spectral density of the process is defined as a Fourier transform (\ref{definition:power_spectral_density}) of the autocovariances $\{ \gamma_h \}_{h \in \mathbb{Z}}$.
Because we consider complex-valued processes, power spectral densities are not generally even functions on $[-\pi, \pi]$.
For the observation $z^{(N)} = (z_1, \cdots, z_N)$ of size $N$ from a Gaussian process with mean zero, let us denote its probability density by $p^{(N)} \bigl( z^{(N)} \bigr)$.
The probability density $p^{(N)} \bigl( z^{(N)} \bigr)$ is explicitly calculated as (\ref{definition:complex_normal_distribution_probability_density_function}) with mean $\mu = 0$ and the following variance-covariance matrix:
\begin{align}
  \Sigma^{(N)}
  := \begin{bmatrix}
    \gamma_0 & \gamma_{-1} & \cdots & \gamma_{-N+1} \\
    \gamma_1 & \gamma_0    & \cdots & \gamma_{-N+2 } \\
    \vdots   & \vdots      & \ddots & \vdots \\
    \gamma_{N-1} & \gamma_{N-2} & \cdots & \gamma_0
  \end{bmatrix}
  \label{definition:variance_covariance_matrix}
  .
\end{align}

As a special case of a strongly stationary Gaussian process with mean zero, we introduce a complex-valued autoregressive moving average (ARMA) process.
A complex-valued ARMA process of degree $(p, q)$ is a Gaussian process that satisfies the relation
\begin{align}
  Z_t = - \sum_{i = 1}^p a_i Z_{t-i} + \varepsilon_{t} + \sum_{i = 1}^q b_i \varepsilon_{t-i}
  \label{definition:ARMA}
\end{align}
for all $t$,
where $a_1, \cdots, a_p, b_1, \cdots, b_q$ are complex-valued coefficients and $\varepsilon_{t}$ is a complex Gaussian white noise with variance $\sigma^2$.
We assume that the polynomials $z^p \bigl( 1 + \sum_{i = 1}^p a_i \, z^{-i} \bigr)$ and $z^q \bigl( 1 + \sum_{i = 1}^q b_i \, z^{-i} \bigr)$ have no common roots in order to ensure identifiability; see Section 3.1 in \cite{brockwell1991time}.
We denote the statistical model of complex-valued stationary ARMA processes by $\mathrm{ARMA}(p, q; \mathbb{C})$ in the present paper.
If $q = 0$, we call the model a complex-valued stationary autoregressive model and denote it by $\mathrm{AR}(p; \mathbb{C})$.
We denote the model of real-valued stationary autoregressive processes of order $p$ by $\mathrm{AR}(p; \mathbb{R})$.
The power spectral density (\ref{definition:power_spectral_density}) of the ARMA model (\ref{definition:ARMA}) is explicitly given by
\begin{align}
  \frac{\sigma^2}{2\pi} \left| \frac{ 1 + \sum_{i = 1}^q b_i \, e^{- i \sqrt{-1} \, \omega} }{ 1 + \sum_{i = 1}^p a_i \, e^{- i \sqrt{-1} \, \omega} } \right|^2
  .
\end{align}

In what follows, we prepare the relevant mathematical notation.
Suppose Gaussian processes are parameterized by complex parameters $\theta \in \Theta \subset \mathbb{C}^p$.
Because we consider complex parameters $\theta \in \Theta \subset \mathbb{C}^p$, we utilize Wirtinger calculus; see Appendix \ref{section:wirtinger_calculus}.
Corresponding to the $i$-th complex coordinate $\theta^i$, there exist Wirtinger derivatives $\pHd{i}$ and $\pAd{i}$.
For simplicity of notation,
we set
\begin{align}
  S_{\alpha_1 \cdots \alpha_a,\, \beta_1 \cdots \beta_b,\, \cdots,\, \gamma_1 \cdots \gamma_c}
  := S^{-1} ( D_{\alpha_1 \cdots \alpha_a} S )
  S^{-1} ( D_{\beta_1 \cdots \beta_b} S )
  \cdots S^{-1} ( D_{\gamma_1 \cdots \gamma_c} S )
\end{align}
for a power spectral density $S = S_{\theta}$ and indices $\alpha_1, \cdots, \alpha_a,\, \beta_1 \cdots \beta_b,\, \cdots,\, \gamma_1 \cdots \gamma_c \in \{1, \cdots, p, \bar{1}, \cdots, \bar{p}\}$,
where $D_{\alpha_1 \cdots \alpha_a} := \pd{\alpha_1} \cdots \pd{\alpha_a} $.
For example, $S_{1\bar{1},2} = S^{-1} (\pHd{1} \pAd{1} S) S^{-1} (\pHd{2} S)$.
We also set
\begin{align}
  M_{\alpha_1 \cdots \alpha_a,\, \beta_1 \cdots \beta_b,\, \cdots,\, \gamma_1 \cdots \gamma_c}
  := \frac{1}{2\pi} \int_{-\pi}^{\pi} S_{\alpha_1 \cdots \alpha_a,\, \beta_1 \cdots \beta_b,\, \cdots,\, \gamma_1 \cdots \gamma_c} (\omega) \,d\omega
  \label{definition:M}
\end{align}
and define the quantities $\gdd{\alpha}{\beta}$, $\Tddd{\alpha}{\beta}{\gamma}$, and $\Guddd{(m)}{\alpha}{\beta}{\gamma}$ as 
\begin{align}
  \gdd{\alpha}{\beta}
  &:= M_{\alpha, \beta}
  = \frac{1}{2 \pi} \int_{-\pi}^{\pi} \left( \pd{\alpha} \log S \right) \left( \pd{\beta} \log S \right) d\omega
  \label{definition:M_g}, \\
  \Tddd{\alpha}{\beta}{\gamma}
  &:= 2 M_{\alpha, \beta, \gamma}
  = \frac{1}{\pi} \int_{-\pi}^{\pi} \left( \pd{\alpha} \log S \right) \left( \pd{\beta} \log S \right) \left( \pd{\gamma} \log S \right) d\omega
  ,\label{definition:M_T}\\
  \Guddd{(m)}{\alpha}{\beta}{\gamma}
  &:= M_{\alpha\beta, \gamma}
  = \frac{1}{2 \pi} \int_{-\pi}^{\pi} \frac{ \pd{\alpha} \pd{\beta} S }{ S } \left( \pd{\gamma} \log S \right) d\omega \label{definition:M_Gm}
\end{align}
for $\alpha, \beta, \gamma \in \{1, \cdots, p, \bar{1}, \cdots, \bar{p} \}$.
The quantities $\Guddd{(m)}{\alpha}{\beta}{\gamma}$ correspond to the coefficients of the mixture connection ($m$-connection) $\nabla^{(m)}$; see Appendix \ref{section:alpha_parallel}.

The complex-valued $2p \times 2p$ matrix $\begin{bmatrix} \gdd{\alpha}{\beta} \end{bmatrix}$ is called the Fisher information matrix,
which naturally induces the metric on the complex parameter space $\Theta$.
The inner product of two functions $N_1 = N_1 \left( \omega \mid \theta \right)$ and $N_2 = N_2 \left( \omega \mid \theta \right)$ defined on $[-\pi, \pi]$ at $\theta \in \Theta$ is defined as
\begin{align}
  \langle N_1, N_2 \rangle_{\theta} := \frac{1}{2\pi} \int_{-\pi}^{\pi} \frac{N_1 \left( \omega \mid \theta \right) }{ S \left( \omega \mid \theta \right) } \frac{N_2 \left( \omega \mid \theta \right) }{ S \left( \omega \mid \theta \right) } \,d\omega
  \,,
\end{align}
and the norm $\| N \|_{\theta}$ of a function $N = N \left( \omega \mid \theta \right)$ at $\theta \in \Theta$ is defined as $\| N \|_{\theta}^2 := \langle N, N \rangle_{\theta}$; see also \cite{komaki1999estimating}.
This form (\ref{definition:M_g}) of the Fisher information matrix was introduced in \cite{whittle1953estimation} for real-valued time-series analysis.
For real-valued processes, the constant $4 \pi$, rather than $2 \pi$, is usually used in the denominator in (\ref{definition:M_g}); see \cite{whittle1953estimation, anderson1977estimation, komaki1999estimating, tanaka2011asymptotic}.
On the other hand, the constant $2 \pi$ is used for signal processing; see \cite{amari1987differential, choi2015kahlerian}.
For complex-valued processes, it is natural to use the constant $2 \pi$ as in (\ref{definition:M_g})
because it yields $\gdd{\alpha}{\beta} = \frac{1}{N} E \bigl[ - \pd{\alpha} \pd{\beta} l^{(N)} \bigr] + O(N^{-1})$, where $l^{(N)}$ denotes the log-likelihood (\ref{definition:log_likelihood}); see Proposition \ref{prop:expectation_loglikehood} in Appendix \ref{proof:expectation_loglikehood}.
Let us denote by $\begin{bmatrix} \guu{\alpha}{\beta} \end{bmatrix}$ the inverse matrix of the Fisher information matrix $\begin{bmatrix} \gdd{\alpha}{\beta} \end{bmatrix}$,
i.e., $\gdd{\alpha}{\gamma} \guu{\gamma}{\beta}  = {\delta_{\alpha}}^{\beta}$ for the Kronecker delta ${\delta_{\alpha}}^{\beta}$.
The prior defined as the square root of the determinant of the $2p \times 2p$ complex-valued matrix $\begin{bmatrix} \gdd{\alpha}{\beta} \end{bmatrix}$ is called the Jeffreys prior and denoted by $\pi_J$ in the present paper.

Asymptotic expansions of Bayesian predictive power spectral densities play an important role in the present paper.
Let us fix a possibly improper prior $\pi$ and assume that $\int_{\Theta} p_{\theta} \bigl( z^{(N)} \bigr) \pi (\theta) \,d\theta$ is finite for any $z^{(N)} \in \mathbb{C}^N$ and that the Bayesian predictive power spectral density (\ref{definition:BPD}) exists for any $\omega \in [-\pi, \pi]$.
The asymptotic expansion of a Bayesian predictive power spectral density (\ref{definition:BPD}) of a complex-valued ARMA process around the maximum likelihood estimator $\hat{\theta} = \hat{\theta}^{(N)} \bigl( z^{(N)} \bigr)$ is
\begin{align}
  \hat{S}_{\pi}^{(N)} (\omega)
  = S \,\bigl( \omega \mid \hat{\theta} \,\bigr)
  + \frac{1}{N} \left( G_{\pi}^{(N)} \,\bigl( \omega \mid \hat{\theta} \,\bigr) + H^{(N)} \,\bigl( \omega \mid \hat{\theta} \,\bigr) \right)
  + O_P \bigl( N^{-\frac{3}{2}} \bigr)
  \label{formula:AEP_S},
\end{align}
where functions $G^{(N)}_{\pi}$ and $H^{(N)}$ represent the parallel and orthogonal parts of the quantity $N \bigl( \hat{S}^{(N)}_{\pi} - S_{\hat{\theta}^{(N)}} \bigr)$, respectively; see Appendix \ref{section:AEP_BPD}.
Functions $G^{(N)}_{\pi}$ and $H^{(N)}$ are explicitly given by
\begin{align}
  G^{(N)}_{\pi} (\omega \mid \theta)
  &:= \guu{\alpha}{\beta} ( \omega \mid \theta ) \left(
    \pd{\alpha} \log \frac{\pi}{\pi_J} (\theta ) + \frac{1}{2} T_{\alpha} ( \omega \mid \theta )
  \right)
  \pd{\beta} S ( \omega \mid \theta )
  \label{definiton:G}, \\
  H^{(N)} (\omega \mid \theta)
  &:= \frac{1}{2} \guu{\alpha}{\beta} ( \omega \mid \theta ) \left(
    \pd{\alpha} \pd{\beta} S ( \omega \mid \theta ) - \Guddu{(m)}{\alpha}{\beta}{\gamma} ( \theta ) \, \pd{\gamma} S ( \omega \mid \theta )
  \right)
  \label{definiton:H},
\end{align}
where $\Guddu{(m)}{\alpha}{\beta}{\gamma} (\theta) := \Guddd{(m)}{\alpha}{\beta}{\delta} (\theta) \, \guu{\delta}{\gamma} (\theta)$ and $T_{\alpha} (\theta) := \Tddd{\alpha}{\beta}{\gamma} (\theta) \, \guu{\beta}{\gamma} (\theta)$.
Note first that $G^{(N)}_{\pi}$ and $H^{(N)}$ are orthogonal in the sense that
$ \langle G^{(N)}_{\pi}, H^{(N)} \rangle_{\theta} = 0 $
for any $\theta \in \Theta$.
Note also that, while the parallel part $G^{(N)}_{\pi}$ may depend on the choice of prior $\pi$, the orthogonal part $H^{(N)}$ is independent of the choice; see \cite{komaki2006shrinkage} for more details.
See Appendix \ref{section:AEP_BPD} for the proof of the expansion in (\ref{formula:AEP_S}).

The Bayesian predictive power spectral density $\hat{S}^{(N)}_{\pi}$ minimizes the Bayes risk (\ref{definition:bayes_risk}) among all the predictive power spectral densities $\hat{S}^{(N)}$ as long as the Bayes risk of $\hat{S}^{(N)}_{\pi}$ is finite; see Appendix \ref{section:kl_divergence}.
Therefore, once we have a prior $\pi$, we can calculate the predictive power spectral density $\hat{S}_{\pi}^{(N)}$ that minimizes the Bayes risk (\ref{definition:bayes_risk}).
The only remaining problem is to find a reasonable prior $\pi$.

\section{K\"ahler parameter spaces for complex-valued autoregressive processes}
\label{section:complex_valued_arma_process}

Let us consider a family $\{ S_{\theta} \}_{\theta \in \Theta}$ of power spectral densities of complex-valued stationary ARMA processes, where $\Theta \subset \mathbb{C}^p$ is a complex parameter space.
If the Fisher information matrix $\begin{bmatrix} \gdd{\alpha}{\beta} \end{bmatrix}$ of the process satisfies the relations
$ \gHdHd{i}{j} = \gAdAd{i}{j} = 0 \,,\, \gHdAd{i}{j} = \gAdHd{j}{i} = \overline{\gHdAd{j}{i}} = \overline{\gAdHd{i}{j}} $
for all $i, j = 1, \cdots, p$ and the relations $ \pHd{i} \gHdAd{j}{k} = \pHd{j} \gHdAd{i}{k} \,,\, \pAd{i} \gHdAd{j}{k} = \pAd{k} \gHdAd{j}{i}$
for all $i, j, k = 1, \cdots, p$, we say that the complex parameter space $\Theta$ is K\"ahler; see also Appendix \ref{section:wirtinger_calculus}.
The K\"ahlerness of the complex parameter space $\Theta$ plays an important role in the construction of priors.

A specific complex parameter space $\Theta \subset \mathbb{C}^{p+q}$ for complex-valued stationary autoregressive moving average processes $\mathrm{ARMA}(p, q; \mathbb{C})$ was shown to be K\"ahler in \cite{choi2015kahlerian}.
We focus on the specific K\"ahler parameter space $\Theta$ for complex-valued stationary autoregressive processes $\mathrm{AR}(p; \mathbb{C})$.
We examine the power spectral densities of $\mathrm{AR}(p; \mathbb{C})$ of the form
\begin{align}
  S(\omega) = \frac{ 1 }{ 2 \pi \left| \prod_{i = 1}^p (1 - \xi^i \, e^{- \sqrt{-1} \, \omega}) \right|^2 }
  \quad, \label{definition:parameter_xi}
\end{align}
where complex parameters $\xi = (\xi^1, \cdots, \xi^p)$ are roots of the polynomial $z^p \bigl( 1 + \sum_{i = 1}^p a_i \, z^{-i} \bigr)$ of the formal variable $z$ and $\sigma^2 = 1$ is assumed.
From the stationarity condition, we assume that $| \xi^i | < 1$ for any $i = 1, \cdots, p$.

We define the parameter space $\tilde{\Theta}_1 \subset \mathbb{C}^p$ as
\begin{align}
  \tilde{\Theta}_1 &:= \left\{ (\xi^1, \cdots, \xi^p) \in \mathbb{C}^p \;\middle|\; \text{$| \xi^i | < 1$ for any $i = 1, \cdots, p$} \right\} \nonumber\\
  &= U \times U \times \cdots \times U,
  \nonumber
\end{align}
where $U$ is the open unit disk in the complex plane $\mathbb{C}$.
In this specific parameterization $\xi = (\xi^1, \cdots, \xi^p)$, the center $0 = (0, \cdots, 0) \in \mathbb{C}^p$ corresponds to the white noise process.
Because we wish to ignore the measure-zero subset, where the denominator of (\ref{definition:parameter_xi}) has multiple roots,
we restrict our attention to the dense subset
\begin{align}
  \Theta_1 := \left\{ (\xi^1, \cdots, \xi^p) \in \tilde{\Theta}_1 \;\middle|\; \text{$\xi^i \neq \xi^j$ for any $i, j = 1, \cdots, p$} \right\}
  \label{definition:parameter_space}
\end{align}
of the original parameter space $\tilde{\Theta}_1$.
The parameter space $\Theta_1$ is a complex manifold of complex dimension $p$ because the set of complex variables $\xi = (\xi^1, \cdots, \xi^p)$ yields a local coordinate of the space,
and the space $\Theta_1$ is open as a topological space with the boundary $\partial \Theta_1$.
In particular, $\Theta_1$ is relatively compact but not compact.

For the specific parameterization $\xi = (\xi^1, \cdots, \xi^p)$ defined in (\ref{definition:parameter_xi}) for $\mathrm{AR} (p, \mathbb{C})$, the Fisher information matrix is explicitly given by
\begin{align}
  \gHdHd{i}{j} = \gAdAd{i}{j} = 0
  \;,\;
  \gHdAd{i}{j} = \gAdHd{j}{i} = \overline{\gHdAd{j}{i}} = \overline{\gAdHd{i}{j}} = \frac{1}{1 - \xi^i \bar{\xi}^j}
\end{align}
for $i, j = 1, \cdots, p$.
Therefore, the complex parameter space $\Theta_1$ is K\"ahler.
This is a very important property of the complex parameter space $\Theta_1$ for analyzing the super-harmonicity of priors.
The $p \times p$ complex-valued matrix $\begin{bmatrix} \gHdAd{i}{j} \end{bmatrix}$ is positive definite if and only if the denominator of (\ref{definition:parameter_xi}) has no multiple roots.
Note that, on the parameter space $\Theta_1$, there exists the inverse $\begin{bmatrix} \gHuAu{i}{j} \end{bmatrix}$ of the Fisher information matrix, which is necessary to define the Laplacian on the parameter space; see Appendix \ref{section:wirtinger_calculus}.

For a K\"ahler parameter space, the Jeffreys prior is the determinant of the $p \times p$ complex-valued Hermitian matrix $\begin{bmatrix} \gHdAd{i}{j} \end{bmatrix}$; see Appendix \ref{section:wirtinger_calculus}.
For the specific parameterization $\xi = (\xi^1, \cdots, \xi^p)$ defined in (\ref{definition:parameter_xi}) for $\mathrm{AR} (p, \mathbb{C})$, the Jeffreys prior $\pi_J$ is explicitly given by
\begin{align}
  \pi_J (\xi) = \frac{\prod_{1 \leq i < j \leq q} | \xi^i - \xi^j |^2 }{\prod_{i = 1}^p \prod_{j = 1}^p \bigl( 1 - \xi^i \bar{\xi}^j \bigr)}
  \label{definition:jeffreys_prior}
  \,;
\end{align}
see also \cite{choi2015kahlerian}.
The Jeffreys prior (\ref{definition:jeffreys_prior}) for $\mathrm{AR} (p, \mathbb{C})$ is continuous in the parameter space $\tilde{\Theta}_1 = U \times \cdots \times U$.
It vanishes if and only if the denominator of (\ref{definition:parameter_xi}) has multiple roots.
Thus, the Jeffreys prior is strictly positive on the parameter space $\Theta_1$.
Moreover, it diverges at the boundary $\partial \tilde{\Theta}_1$ of the parameter space $\tilde{\Theta}_1$ and defines an improper prior on $\tilde{\Theta}_1$.

\section{Main Theorem}
\label{section:main_theorem}

Let us consider a family $\{ S_{\theta} \}_{\theta \in \Theta}$ of power spectral densities of complex-valued stationary processes, where $\Theta \subset \mathbb{C}^p$ is a complex parameter space.
Our objective is to construct a predictive power spectral density $\hat{S}^{(N)}$,
the risk $R \bigl( \hat{S}^{(N)} \bigr)$ of which is kept as small as possible.
We say that a predictive power spectral density $\hat{S}^{(N)}_1$ dominates a predictive power spectral density $\hat{S}^{(N)}_2$ if $R \bigl( \hat{S}^{(N)}_1 \mid \theta \bigr) \leq R \bigl( \hat{S}^{(N)}_2 \mid \theta \bigr)$ for any $\theta \in \Theta$ and the strict inequality holds for some $\theta$.

Suppose that the parameter space $\Theta$ is K\"ahler and there exists a positive continuous eigenfunction $\phi$ of the Laplacian (Laplace--Beltrami operator; see (\ref{formula:laplacian}) in Appendix \ref{section:wirtinger_calculus}) $\Delta := 2 \gHuAu{i}{j} \pHd{i} \pAd{j}$ with a negative eigenvalue $-K$ globally defined on $\Theta$.
We define a family of priors $\{ \pi^{(\kappa)} \}_{\kappa \in \mathbb{R}}$, called $\kappa$-priors, as $\pi^{(\kappa)} := \phi^{-\kappa + 1} \pi_J$, where $\pi_J$ denotes the Jeffreys prior.
We state that, with a suitable choice of $\kappa \in \mathbb{R}$, the Bayesian predictive power spectral density $\hat{S}^{(N)}_{\pi^{(\kappa)}}$ based on the proposed prior $\pi^{(\kappa)}$ asymptotically dominates the Bayesian predictive power spectral density $\hat{S}^{(N)}_{\pi_J}$ based on the Jeffreys prior $\pi_J$.
We have proved the following theorem, which is particularly for $\kappa$-priors for complex-valued ARMA processes.
\begin{thm}[Main Theorem]
  \label{theorem:main}
  Let $\{ S_{\theta} \}_{\theta \in \Theta}$ be a family of power spectral densities parameterized by a K\"ahler parameter space $\Theta \subset \mathbb{C}^p$.
  We suppose that the inverse $\begin{bmatrix} \gHuAu{i}{j} \end{bmatrix}$ of the Fisher information matrix $\begin{bmatrix} \gHdAd{i}{j} \end{bmatrix}$ exists on $\Theta$ and, hence, the Laplacian $\Delta$ is properly defined on $\Theta$.
  We also suppose that there exists a positive continuous eigenfunction $\phi$ of the Laplacian $\Delta$ with a negative eigenvalue $-K$ globally defined on $\Theta$.
  Let $\pi_1 := \pi^{(\kappa_1)}$ and $\pi_2 := \pi^{(\kappa_2)}$ be two $\kappa$-priors for $\kappa_1, \kappa_2 \in \mathbb{R}$,
  and assume that Bayesian predictive power spectral densities $\hat{S}^{(N)}_{\pi^{(\kappa)}}$ exist for $\kappa = \kappa_1, \kappa_2$.
  Further assume that there exist the asymptotic expansions (\ref{formula:AEP_S}) of the Bayesian predictive power spectral densities around the maximum likelihood estimator.
  Then, we have
  \begin{align*}
    R \bigl( \hat{S}^{(N)}_{\pi_1} \mid \theta \bigr) - R \bigl( \hat{S}^{(N)}_{\pi_2} \mid \theta \bigr)
    = \frac{(\kappa_1 - \kappa_2) K + (\kappa_1 - \kappa_2) (\kappa_1 + \kappa_2) \gHuAu{i}{j} \bigl( \pHd{i} \log \phi \bigr) \bigl( \pAd{j} \log \phi \bigr)}{N^2}
    + O \bigl( N^{-\frac{5}{2}} \bigr)
  \end{align*}
  for $\theta \in \Theta$.
\end{thm}

The proof of the Main Theorem is largely aided by the form (\ref{formula:AEP_S}) of the asymptotic expansion of the Bayesian predictive power spectral density $\hat{S}^{(N)}_{\pi}$ around the maximum likelihood estimator $\hat{\theta}^{(N)}$.
The comprehensive proof of the Main Theorem is given in Appendix \ref{proof:main}.
Note that the specific construction of the K\"ahler parameter space (\ref{definition:parameter_space}) for $\mathrm{AR} (p; \mathbb{C})$ is unrelated to this theorem.
This theorem is a very general theorem that always holds as long as the parameter space $\Theta$ is K\"ahler.
The importance of the K\"ahlerness of parameter spaces in statistics and the generalization of the Main Theorem for the i.i.d.\ case are discussed in Section \ref{section:generalization}.
The metric $\gHdAd{i}{j}$ and its inverse $\gHuAu{i}{j}$ of a K\"ahler manifold are explained in Appendix \ref{section:wirtinger_calculus}.

Setting $\pi_1 := \pi^{(+1)} = \pi_J$ in Theorem \ref{theorem:main}, we can easily observe that if $-1 \leq \kappa < 1$, then $\hat{S}^{(N)}_{\pi^{(\kappa)}}$ asymptotically dominates $\hat{S}^{(N)}_{\pi_J}$.

\begin{corollary}
  \label{corollary:main}
  Let $\pi := \pi^{(\kappa)}$ be a $\kappa$-prior on a K\"ahler parameter space $\Theta \subset \mathbb{C}^p$.
  We have
  \begin{align*}
    R \bigl( \hat{S}^{(N)}_{\pi_J} \mid \theta \bigr) - R \bigl( \hat{S}^{(N)}_{\pi^{(\kappa)}} \mid \theta \bigr)
    = \frac{(1 - \kappa) K + (1 - \kappa^2) \gHuAu{i}{j} \bigl( \pHd{i} \log \phi \bigr) \bigl( \pAd{j} \log \phi \bigr)}{N^2}
    + O \bigl( N^{-\frac{5}{2}} \bigr)
    .
  \end{align*}
  Therefore, if $-1 \leq \kappa < 1$, then $\hat{S}^{(N)}_{\pi^{(\kappa)}}$ asymptotically dominates $\hat{S}^{(N)}_{\pi_J}$.
\end{corollary}

Recall that $\gHuAu{i}{j} \bigl( \pHd{i} \log \phi \bigr) \bigl( \pAd{j} \log \phi \bigr) \geq 0$ because the Hermitian matrix $\begin{bmatrix} \gHuAu{i}{j} \end{bmatrix}$ is positive definite.
Therefore, to maximize the worst case of the risk improvement, we propose the $\kappa$-prior for $\kappa = -1$.
When $\kappa = -1$, the Bayesian predictive power spectral density $\hat{S}^{(N)}_{\pi^{(-1)}}$ achieves constant risk improvement:
$
R \bigl( \hat{S}^{(N)}_{\pi_J} \mid \theta \bigr) - R \bigl( \hat{S}^{(N)}_{\pi^{(-1)}} \mid \theta \bigr) 
= 2 K / N^2
+ O \bigl( N^{-\frac{5}{2}} \bigr)
$.
Equation (\ref{formula:main_ar}) is a special case of Corollary \ref{corollary:main} for $\mathrm{AR} (p; \mathbb{C})$ and $\kappa = -1$,
where $K = p(p+1)$.
Section \ref{section:super_harmonicity_prior} discusses the existence of the positive continuous eigenfunction $\phi$ with a negative eigenvalue $-K = -p(p+1)$ on the parameter space $\Theta_1$ for $\mathrm{AR} (p; \mathbb{C})$.

\section{Super-harmonic priors on $\mathrm{AR} (p; \mathbb{C})$}
\label{section:super_harmonicity_prior}

In this section, we prove the existence of the positive continuous eigenfunction $\phi$ of the Laplacian $\Delta$ with eigenvalue $-K = -p(p+1)$ for $\mathrm{AR} (p; \mathbb{C})$.
Furthermore, for the $\mathrm{AR} (p; \mathbb{C})$ model, we show that the Bayesian predictive power spectral density $\hat{S}^{(N)}_{\pi^{(-1)}}$ based on the $(-1)$-prior $\pi^{(-1)}$ asymptotically dominates the estimative power spectral density $S_{\hat{\theta}^{(N)}}$ with the maximum likelihood estimator $\hat{\theta}^{(N)}$.
This is another reason why we propose the $(-1)$-prior $\pi^{(-1)}$ for the case of $\mathrm{AR} (p; \mathbb{C})$.

The eigenfunction $\phi$ for $\mathrm{AR} (p; \mathbb{C})$ is defined as
\begin{align}
  \phi (\xi) := \prod_{i=1}^p \prod_{j=1}^p \bigl( 1 - \xi^i \bar{\xi}^j \bigr)
  \label{definition:phi}
\end{align}
for $\xi = (\xi^1, \cdots, \xi^p) \in \tilde{\Theta}_1 = U \times \cdots \times U$.
The function $\phi$ is the inverse of the determinant $\bigl| \Sigma^{(N)} \bigr|$ of the variance-covariance matrix $\Sigma^{(N)}$ of size $N \geq p$ for $\mathrm{AR} (p; \mathbb{C})$; see Appendix \ref{section:existenece_BPSD}.
The function $\phi$ is a real-valued continuous function defined globally on the parameter space $\tilde{\Theta}_1$.
Moreover, it is positive on $\tilde{\Theta}_1$ and is $0$ at the boundary $\partial \tilde{\Theta}_1$ of $\tilde{\Theta}_1$.
Note also that the function $\phi$ has its maximum at the white noise process.

The $\kappa$-prior $\pi^{(\kappa)}$ for $\mathrm{AR} (p; \mathbb{C})$ is
\begin{align}
  \pi^{(\kappa)} := \phi^{-\kappa+1} \pi_J
  = \left( \prod_{i=1}^p \prod_{j=1}^p \bigl( 1 - \xi^i \bar{\xi}^j \bigr) \right)^{-\kappa} 
  \left( \prod_{1 \leq i < j \leq q} | \xi^i - \xi^j |^2 \right)
  ,
  \label{definition:alpha_prior}
\end{align}
where $\pi_J$ is the Jeffreys prior (\ref{definition:jeffreys_prior}) for $\mathrm{AR} (p; \mathbb{C})$.
The $\kappa$-prior $\pi^{(\kappa)}$ for $\mathrm{AR} (p; \mathbb{C})$ is proper if $\kappa < 1$ and improper if $\kappa \geq 1$ on the parameter space $\Tilde{\Theta}_1$; see Appendix \ref{section:existenece_BPSD}.
In particular, the Jeffreys prior $\pi_J = \pi^{(+1)}$ is improper on the parameter space $\Tilde{\Theta}_1$.
The Bayesian predictive power spectral densities $\hat{S}^{(N)}_{\pi}$ for $\mathrm{AR} (p; \mathbb{C})$ based on the $\kappa$-prior $\pi^{(\kappa)}$ exists if $\kappa < 2$ and $N \geq p$; see Appendix \ref{section:existenece_BPSD}.

Before proving $\Delta \phi = -p(p+1) \phi$, we introduce a useful lemma.
\begin{lemma}
  \label{lemma:phi}
  We have
  \begin{align}
    &\pHd{i} \log \phi = \sum_{j=1}^p \frac{-\bar{\xi}^j}{1 - \xi^i \bar{\xi^j}} = - \gHdAd{i}{j} \bar{\xi}^j, \\
    &\pAd{j} \log \phi = \sum_{i=1}^p \frac{-\xi^i}{1 - \xi^i \bar{\xi^j}} = - \gHdAd{i}{j} \xi^i, 
  \end{align}
  and
  \begin{align}
    \xi^i \pHd{i} \log \pi_J = \bar{\xi}^j \pAd{j} \log \pi_J = \frac{1}{2} p (p-1) + \xi^i \bar{\xi}^j \gHdAd{i}{j} .
  \end{align}
\end{lemma}

\begin{IEEEproof}
  To prove the third equation, use the identity
  $ \sum_{i=1}^p \sum_{j=1, j \neq i}^p \frac{\xi^i}{\xi^i - \xi^j} = \frac{1}{2} p (p - 1)$ .
\end{IEEEproof}

Using Lemma \ref{lemma:phi}, we see that $\phi$ is, in fact, an eigenfunction of the Laplacian with eigenvalue $-K = -p(p+1)$.

\begin{prop}
  \begin{align}
    \Delta \phi 
    &= -p (p+1) \phi
    .
  \end{align}
\end{prop}

\begin{IEEEproof}
  Because the parameter space is K\"ahler, we can use the formula (\ref{formula:laplacian}) for its definition of the Laplacian.
  Direct computation shows
  \begin{align*}
    \Delta \phi 
    &= 2 \gHuAu{i}{j} \pHd{i} \pAd{j} \phi = 2 \gHuAu{i}{j} \pHd{i} \left( - \gHdAd{k}{j} \xi^k \phi \right)\\
    &= 2 \gHuAu{i}{j} \left(
      - (\pHd{i} \gHdAd{k}{j} ) \xi^k \phi - \gHdAd{k}{j} ( \pHd{i} \xi^k) \phi - \gHdAd{k}{j} \xi^k (\pHd{i} \phi))
    \right) \\
    &= 2 \gHuAu{i}{j} \left(
      - (\pHd{k} \gHdAd{i}{j} ) \xi^k \phi - \gHdAd{i}{j} \phi + \gHdAd{k}{j} \gHdAd{i}{l} \xi^k \bar{\xi}^l \phi)
    \right) \\
    &= -2 \, \xi^k (\pHd{k} \log \pi_J) \phi - 2 \gHuAu{i}{j} \gHdAd{i}{j} \phi + 2 \gHdAd{k}{l} \xi^k \bar{\xi}^l \phi \\
    &= - p (p - 1) \phi - 2 p \phi = - p (p + 1) \phi 
    ,
  \end{align*}
  where we have used the K\"ahlerness (\ref{definition:kaehler_manifold}), the Jacobi formula (\ref{formula:jacobi}), and $\gHuAu{i}{j} \gHdAd{i}{j} = p$.
\end{IEEEproof}

As stated in Corollary \ref{corollary:main}, $\hat{S}^{(N)}_{\pi^{(\kappa)}}$ asymptotically dominates $\hat{S}^{(N)}_{\pi^{(+1)}} = \hat{S}^{(N)}_{\pi_J}$ if $-1 \leq \kappa < 1$.
For $\mathrm{AR} (2; \mathbb{C}$), the specific prior $\psi := (1 - \xi^1 \bar{\xi}^2) (1 - \xi^2 \bar{\xi}^1)(1 - |\xi^1|^2)(1 - |\xi^2|^2)$ is introduced as a super-harmonic prior in \cite{choi2015kahlerian}.
This prior $\psi$ is the special case of Corollary \ref{corollary:main} for $(p, \kappa) = (2, 0)$.
For $\mathrm{AR} (p; \mathbb{R})$ with $p \geq 2$, a similar but slightly different prior is presented in \cite{tanaka2018superharmonic}.
This prior corresponds to the $\kappa = 0$ case for a positive eigenfunction of the Laplacian with eigenvalue $-K = -p(p-1)$.

We show that the Bayesian predictive power spectral density $\hat{S}^{(N)}_{\pi^{(-1)}}$ asymptotically dominates the estimative power spectral density $S_{\hat{\theta}^{(N)}}$ with the maximum likelihood estimator $\hat{\theta}^{(N)}$.
Let us fix the true parameter $\theta_0 \in \Theta$ and denote the maximum likelihood estimator by $\hat{\theta}^{(N)} := \hat{\theta}^{(N)} \bigl( z^{(N)} \bigr)$ for the observation $z^{(N)} \in \mathbb{C}^N$.
According to \cite{komaki1999estimating},
if we fortunately find a prior $\pi$ such that $G^{(N)}_{\pi} = 0$,
then we have
\begin{align}
R \bigl( S_{\hat{\theta}^{(N)}} \mid \theta_0 \bigr) - R \bigl( \hat{S}^{(N)}_{\pi} \mid \theta_0 \bigr) 
= \frac{1}{2 N^2} \| H^{(N)} \|^2_{\theta_0} + O (N^{-\frac{5}{2}})
.
\end{align}
For the specific parametrization $\xi \in \Theta_1 \subset \mathbb{C}^p$ for $\mathrm{AR} (p; \mathbb{C})$ defined in (\ref{definition:parameter_xi}),
direct computation shows
\begin{align}
  G^{(N)}_{\pi^{(\kappa)}} &= 2 \gHuAu{i}{j} \left(
    \pHd{i} \log \frac{\pi^{(\kappa)}}{\pi_J} + \frac{1}{2} T_{i}
  \right)
  \pAd{j} S
  = - 2 (\kappa + 1) \bar{\xi}^j \pAd{j} S
\end{align}
for the $\kappa$-prior $\pi^{(\kappa)} := \phi^{-\kappa+1} \pi_J$.
Thus, if $\kappa = -1$, then $G^{(N)}_{\pi} = 0$.
Therefore, the Bayesian predictive power spectral density $\hat{S}^{(N)}_{\pi^{(-1)}}$ asymptotically dominates the estimative power spectral density $S_{\hat{\theta}^{(N)}}$ with the maximum likelihood estimator $\hat{\theta}^{(N)}$.

\section{Generalization of the Main Theorem}
\label{section:generalization}

Although the present study mainly focuses on the complex Gaussian process, the Main Theorem is valid for the i.i.d.\ case as long as the complex parameter space is K\"ahler.

K\"ahlerness is a complexified concept of the exponential family.
Consider a family of probability density functions $\{ p_{\theta} \}_{\theta \in \Theta}$ of the exponential family parameterized by real parameters $\theta = (\theta^1, \cdots, \theta^p) \in \Theta \subset \mathbb{R}^p$.
We may assume that the probability density function is of the form $p_{\theta} (x) = \exp ( \theta^i x_i - \Psi(\theta))$.
We know that the Fisher information matrix is given by $\gdd{i}{j} = \pd{i} \pd{j} \Psi$.
The K\"ahler parameter space $\Theta$ is the generalization of the exponential family in the sense that there exists, at least locally, a function $\mathcal{K}$ on $\Theta$, called a K\"ahler potential, such that $\gHdAd{i}{j} = \pHd{i} \pAd{j} \mathcal{K}$.
If there exists a K\"ahler potential $\mathcal{K}$ on $\Theta$, then it is easy to see that the definition of K\"ahlerness (\ref{definition:kaehler_manifold}) holds.
The converse is also true; see \cite{choi2015kahlerian, nakahara2003geometry}.

We define the risk of a Bayesian predictive distribution for the i.i.d\ case.
Consider a family of probability density functions $\{ p_{\theta} \}_{\theta \in \Theta}$ parameterized by complex parameters $\theta = (\theta^1, \cdots, \theta^p) \in \Theta \subset \mathbb{C}^p$, where $\Theta$ is K\"ahler.
The sample space $\mathcal{Z}$ of this model may be any subset of $\mathbb{R}^r$ or $\mathbb{C}^r$.
Let $\pi$ be a possibly improper prior for this model.
Consider an i.i.d.\ sample $z^{(N)} = (z_1, \cdots, z_N) \in \mathcal{Z}^N$ of size $N$ from the distribution at $\theta \in \Theta$.
The predictive distribution $\hat{p}^{(N)}_{\pi} (z) := \int_{\Theta} p_{\theta} (z) \, \pi \bigl( \theta \mid z^{(N)} \bigr) \,d\theta$ for $z \in \mathcal{Z}$ is called the Bayesian predictive distribution based on the prior $\pi$.
The risk of $\hat{p}^{(N)}_{\pi}$ is defined as
\begin{align*}
  R \bigl( \hat{p}^{(N)}_{\pi} \mid \theta \bigr)
  &:= E_{\theta} \left[ \KL{p_{\theta}}{\hat{p}^{(N)}_{\pi}} \right] \\
  &= \int_{\mathcal{Z}^N} \KL{p_{\theta}}{\hat{p}^{(N)}_{\pi}} dP^{(N)}_{\theta} \bigl( z^{(N)} \bigr) \\
  &= \int_{\mathcal{Z}^N}
  \left( \int_{\mathcal{Z}} p_{\theta} (z) \log \frac{p_{\theta} (z)}{\hat{p}^{(N)}_{\pi} (z)} dz \right)
  dP^{(N)}_{\theta} \bigl( z^{(N)} \bigr)
  ,
\end{align*}
where $dP^{(N)}_{\theta} \bigl( z^{(N)} \bigr) := \bigl( \prod_{t = 1}^{N} p_{\theta} (z_t) \,dz_t \bigr)$.

We prepare the relevant mathematical notation.
We set
\begin{align}
  K_{\alpha_1 \cdots \alpha_a,\, \beta_1 \cdots \beta_b,\, \cdots,\, \gamma_1 \cdots \gamma_c}
  := \int_{\mathcal{Z}} \bigl( D_{\alpha_1 \cdots \alpha_a} \, \log p (z) \bigr) \, \bigl( D_{\beta_1 \cdots \beta_b} \, \log p (z) \bigr)
  \cdots \bigl( D_{\gamma_1 \cdots \gamma_c} \, \log p (z) \bigr) \,P (dz)
  \label{definition:K}
\end{align}
for a probability density function $p = p_{\theta}$ and indices $\alpha_1, \cdots, \alpha_a,\, \beta_1 \cdots \beta_b,\, \cdots,\, \gamma_1 \cdots \gamma_c \in \{1, \cdots, p, \bar{1}, \cdots, \bar{p}\}$,
where $P (dz) := p (z) \, dz$.
We define the quantities $\gdd{\alpha}{\beta}$, $\Tddd{\alpha}{\beta}{\gamma}$, and $\Guddd{(m)}{\alpha}{\beta}{\gamma}$ as 
\begin{align}
  \gdd{\alpha}{\beta}
  &:= K_{\alpha, \beta} = E \left[ (\pd{\alpha} \log p) (\pd{\beta} \log p) \right], \label{definition:K_g}\\
  \Tddd{\alpha}{\beta}{\gamma}
  &:= K_{\alpha, \beta, \gamma} = E \left[ (\pd{\alpha} \log p) (\pd{\beta} \log p) (\pd{\gamma} \log p) \right], \label{definition:K_T}\\
  \Guddd{(m)}{\alpha}{\beta}{\gamma}
  &:= K_{\alpha\beta, \gamma} + K_{\alpha, \beta, \gamma} = E \left[ (\pd{\alpha} \pd{\beta} \, p) (\pd{\gamma} \, p) \right] \label{definition:K_Gm}
  ,
\end{align}
respectively, for $\alpha, \beta, \gamma \in \{1, \cdots, p, \bar{1}, \cdots, \bar{p} \}$.
The relation to the former definitions of the quantities $\gdd{\alpha}{\beta}$, $\Tddd{\alpha}{\beta}{\gamma}$, and $\Guddd{(m)}{\alpha}{\beta}{\gamma}$, namely (\ref{definition:M_g}), (\ref{definition:M_T}), and (\ref{definition:M_Gm}), is given in Proposition \ref{prop:expectation_loglikehood} in Appendix \ref{proof:expectation_loglikehood}; they coincide with each other up to $O(1)$ terms for complex-valued stationary ARMA processes.

We have the same form of the asymptotic expansion of $\hat{p}^{(N)}_{\pi}$ as (\ref{formula:AEP_S}) for the i.i.d.\ case; see \cite{komaki1996asymptotic}.
Consider a positive continuous eigenfunction $\phi$ of the Laplacian $\Delta$ with a negative eigenvalue $-K$, i.e., $\Delta \phi = -K \phi < 0$.
Then, we can construct the $\kappa$-prior by $\pi^{(\kappa)} := \phi^{-\kappa+1} \pi_J$ for $\kappa \in \mathbb{R}$, where $\pi_J$ is the Jeffreys prior of this model.
As the proof of the Main Theorem (Proposition \ref{theorem:main}) only depends on the form (\ref{formula:AEP_S}) of the asymptotic expansion of the Bayesian predictive power spectral density $\hat{S}^{(N)}_{\pi}$,
Theorem \ref{theorem:main} also holds for the risk difference $R \bigl( \hat{p}^{(N)}_{\pi_1} \mid \theta \bigr) - R \bigl( \hat{p}^{(N)}_{\pi_2} \mid \theta \bigr)$ for the i.i.d.\ case.
Therefore, our proposal is to use the prior $\pi^{(-1)} := \phi^{-2} \pi_J$, where $\phi$ is an eigenfunction of the Laplacian with the smallest negative eigenvalue.

The construction $\pi^{(\kappa)} := \phi^{-\kappa+1} \pi_J$ of a family of $\kappa$-priors is related to other types of objective priors.
For example, if $\pHd{i} \phi$ is proportional to $\THd{i} := \THdHdAd{i}{k}{j} \gHuAu{k}{j} + \THdAdHd{i}{j}{k} \gHuAu{k}{j}$ on the parameter space, the family $\{ \pi^{(\kappa)} \}_{\kappa \in \mathbb{R}}$ of $\kappa$-priors is, in fact, a family of $\alpha$-parallel priors; see Appendix \ref{section:alpha_parallel}.
In general, $\alpha$-parallel priors do not always exist.
Therefore, the existence of a family of $\alpha$-parallel priors suggests some statistical property in the statistical model; see \cite{takeuchi2005spl}.

\section{Numerical experiments for risk differences}
\label{section:numerical_experiments}

We consider the $\mathrm{AR} (1; \mathbb{C})$ case of $z^t = \xi z^{t-1} + \varepsilon$ in this section.
The parameter space is the open unit disk $U = \{ \xi \in \mathbb{C} \mid |\xi| < 1 \}$.
The $\kappa$-prior for $\mathrm{AR} (1; \mathbb{C})$ is
$ \pi^{(\kappa)} (\xi) := (1 - |\xi|^2)^{-\kappa} $,
where $\xi \in U$.
Recall that $\kappa = +1$ corresponds to the improper Jeffreys prior $\pi_J = (1 - |\xi|^2)^{-1}$, which is mentioned as a reference prior in \cite{berger1994noninformative} for $\mathrm{AR} (1, \mathbb{R})$.
On the other hand, $\kappa = -1$ corresponds to the proposed proper prior $\pi^{(-1)} = (1 - |\xi|^2)$, 
which is the inverse of the Jeffreys prior, and is also mentioned in \cite{zellner1977maximal} for the $\mathrm{AR} (1; \mathbb{R})$ case.
In fact, the inverse of the Jeffreys prior for $\mathrm{AR} (1; \mathbb{R})$ is a maximal data information prior (MDIP) for $\mathrm{AR} (1; \mathbb{R})$; see \cite{zellner1977maximal}.
Note that for $p \geq 2$, the proposed proper prior $\pi^{(-1)}$ is not the inverse of the Jeffreys prior $\pi^{(+1)}$ for $\mathrm{AR} (p; \mathbb{C})$.

\begin{figure}[htbp]
  \centering
  \includegraphics[width=0.5\columnwidth]{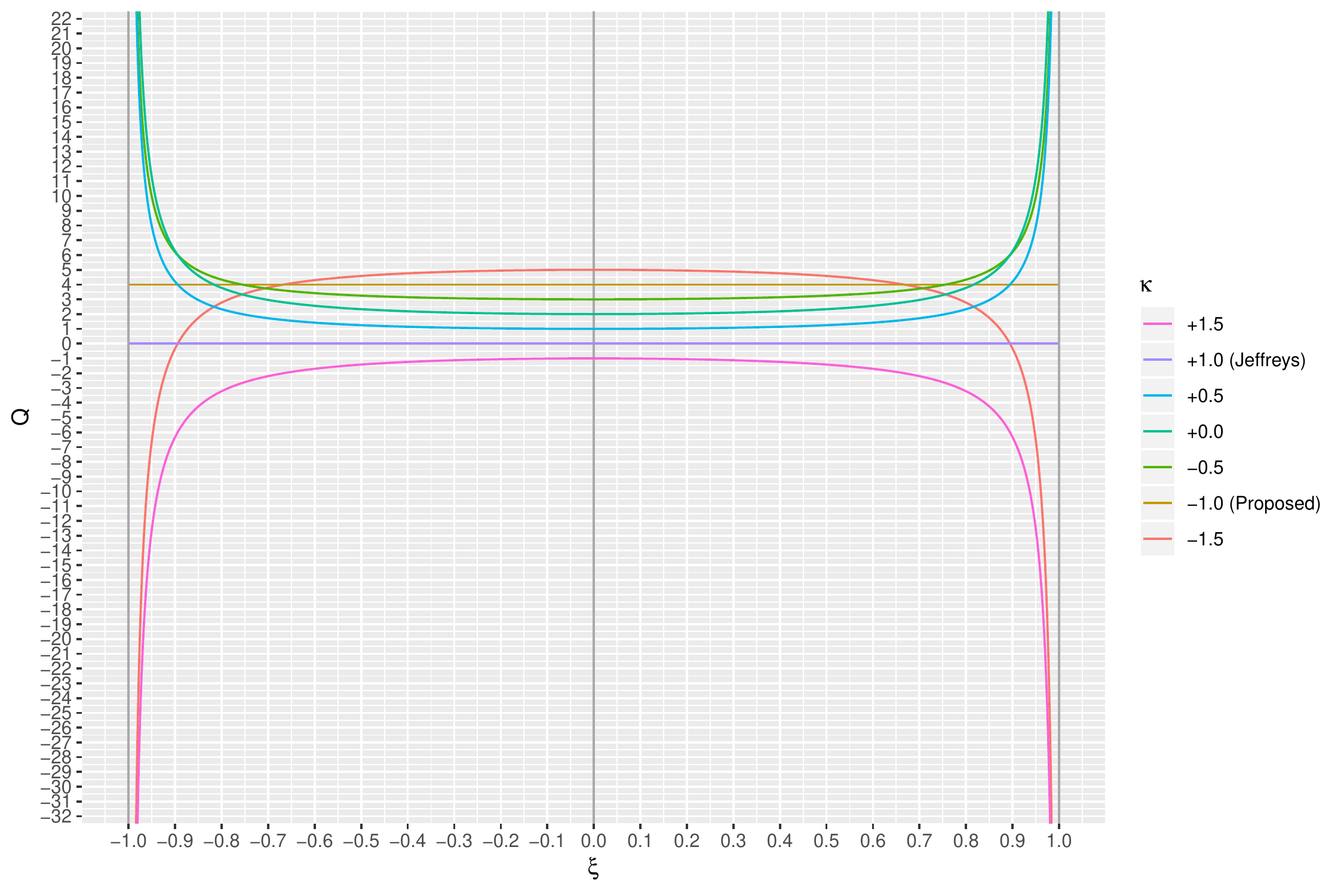}
  \caption{Expected pointwise limit $Q^{(\kappa)} (\xi)$ for the risk difference $Z^{(\kappa)} (\xi)$ for $\mathrm{AR} (1; \mathbb{C})$, where $\xi \in (-1, 1) \subset U$.}
  \label{figure:KL_DIFF_LIMIT}
\end{figure}

Corollary \ref{corollary:main} for $\mathrm{AR} (1; \mathbb{C})$ now becomes
\begin{align}
  R \bigl( \hat{S}^{(N)}_{\pi_J} \mid \xi \bigr) - R \bigl( \hat{S}^{(N)}_{\pi^{(\kappa)}} \mid \xi \bigr)
  = \frac{1}{N^2} Q^{(\kappa)} (\xi) + O \bigl( N^{-\frac{5}{2}} \bigr)
  ,
\end{align}
where
$ Q^{(\kappa)} (\xi) := 2(1 - \kappa) + (1 - \kappa^2) \frac{|\xi|^2}{1 - |\xi|^2} $
is the expected pointwise limit of the normalized risk difference
$ Z^{(\kappa)} (\xi) := N^2 \left( R \bigl( \hat{S}^{(N)}_{\pi_J} \mid \xi \bigr) - R \bigl( \hat{S}^{(N)}_{\pi^{(\kappa)}} \mid \xi \bigr) \right)$.
Figure \ref{figure:KL_DIFF_LIMIT} shows the function $Q^{(\kappa)} (\xi)$.
As stated in Corollary \ref{corollary:main}, if $-1 \leq \kappa \leq 1$, then $Q^{(\kappa)} (\xi) \geq 0$.
In particular, if $\kappa = -1$, then the risk difference $Z^{(\kappa)} (\xi)$ asymptotically achieves the constant $Q^{(-1)} (\xi) = 2K = 4$.

\begin{figure}[htbp]
  \centering
  \includegraphics[width=0.5\columnwidth]{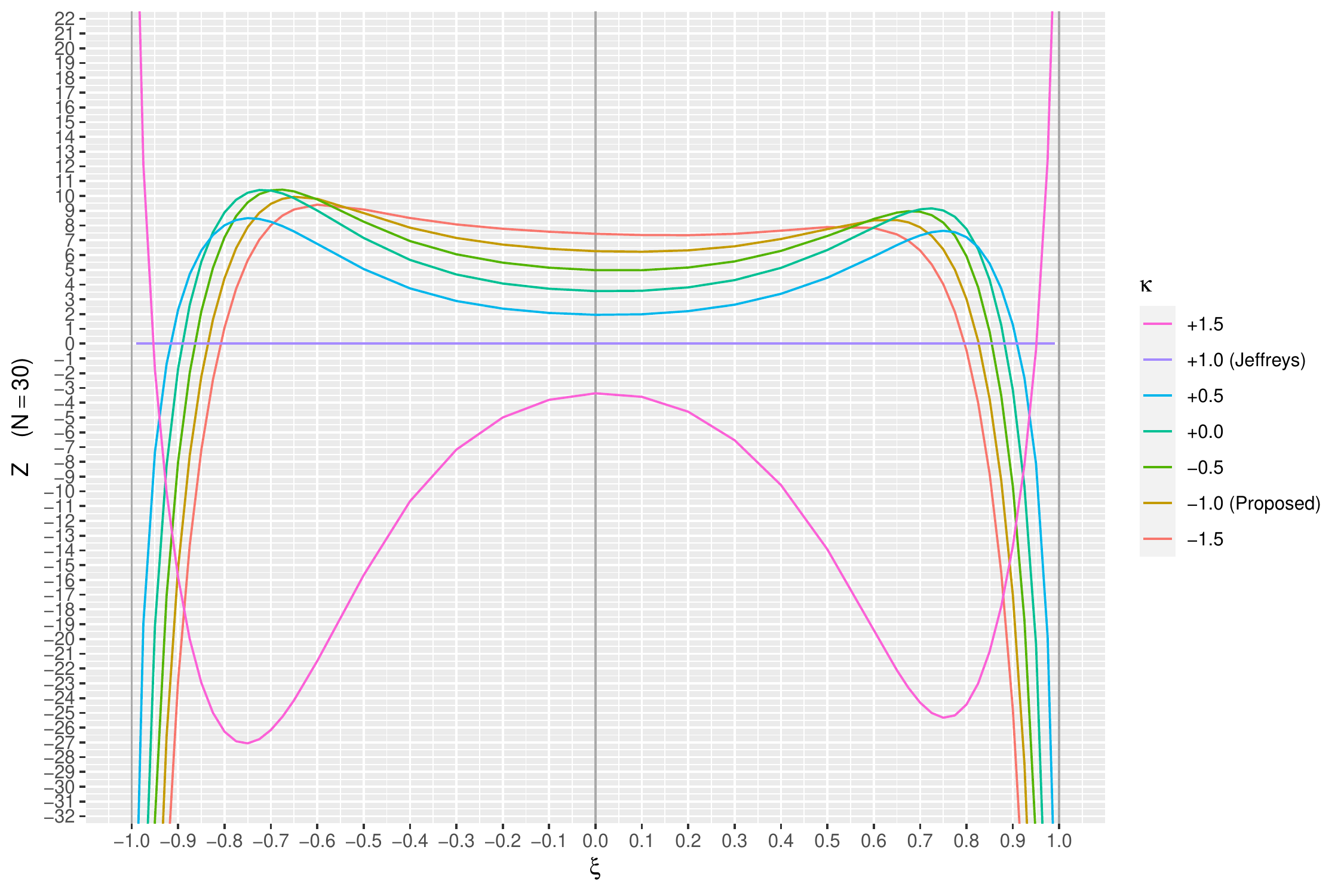}
  \caption{Numerical experiment for the risk difference $Z^{(\kappa)} (\xi)$ for $\mathrm{AR} (1; \mathbb{C})$ with $N = 30$, where $\xi \in (-1, 1) \subset U$. The proposed $\hat{S}^{(N)}_{\pi^{(-1)}}$ dominates the baseline $\hat{S}^{(N)}_{\pi_J}$ for $\xi \in [-0.825, +0.825]$.}
  \label{figure:KL_DIFF_N30}
\end{figure}

\begin{figure}[htbp]
  \centering
  \includegraphics[width=0.5\columnwidth]{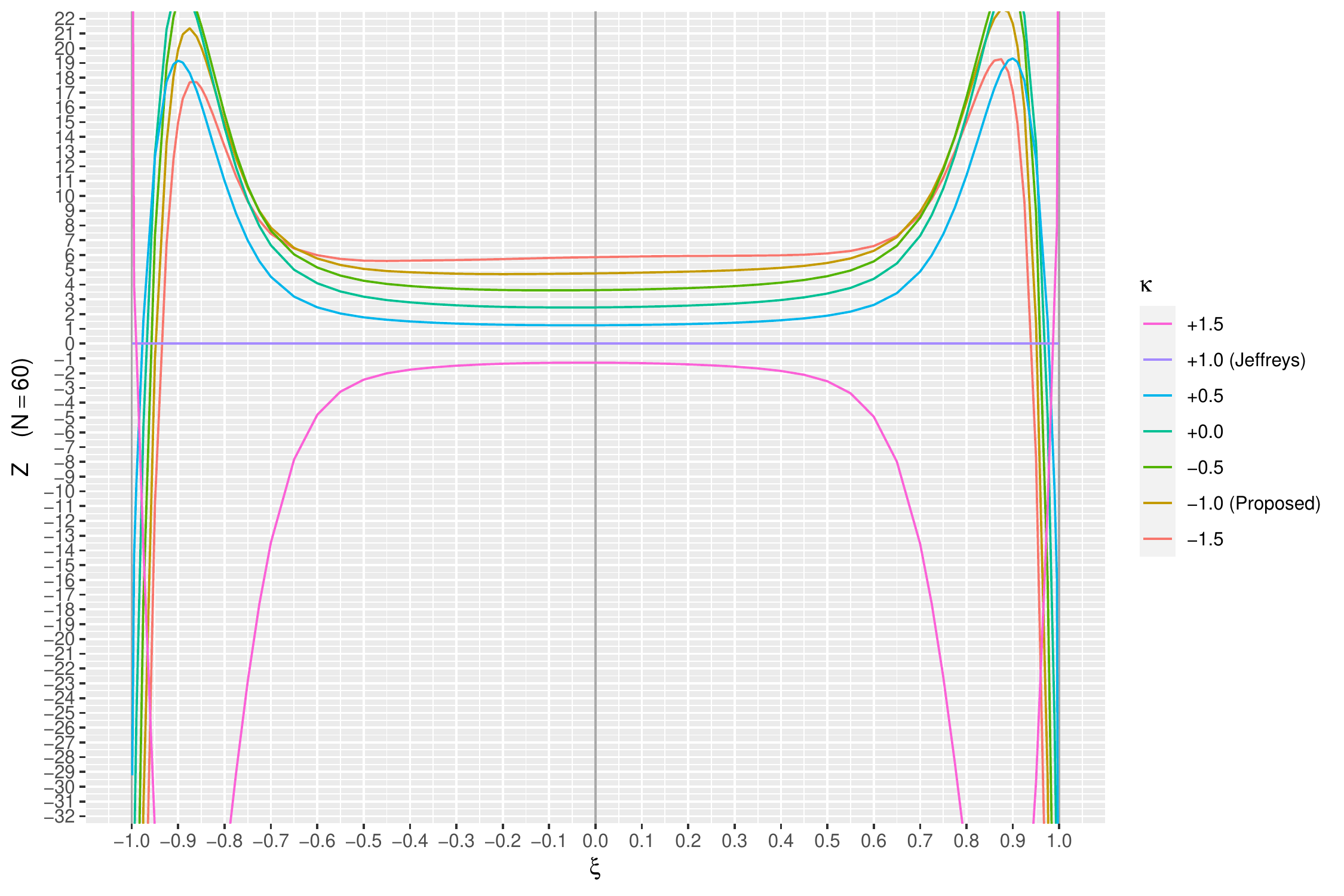}
  \caption{Numerical experiment for the risk difference $Z^{(\kappa)} (\xi)$ for $\mathrm{AR} (1; \mathbb{C})$ with $N = 60$, where $\xi \in (-1, 1) \subset U$. The proposed $\hat{S}^{(N)}_{\pi^{(-1)}}$ dominates the baseline $\hat{S}^{(N)}_{\pi_J}$ for $\xi \in [-0.925, +0.925]$.}
  \label{figure:KL_DIFF_N60}
\end{figure}

\begin{figure}[htbp]
  \centering
  \includegraphics[width=0.5\columnwidth]{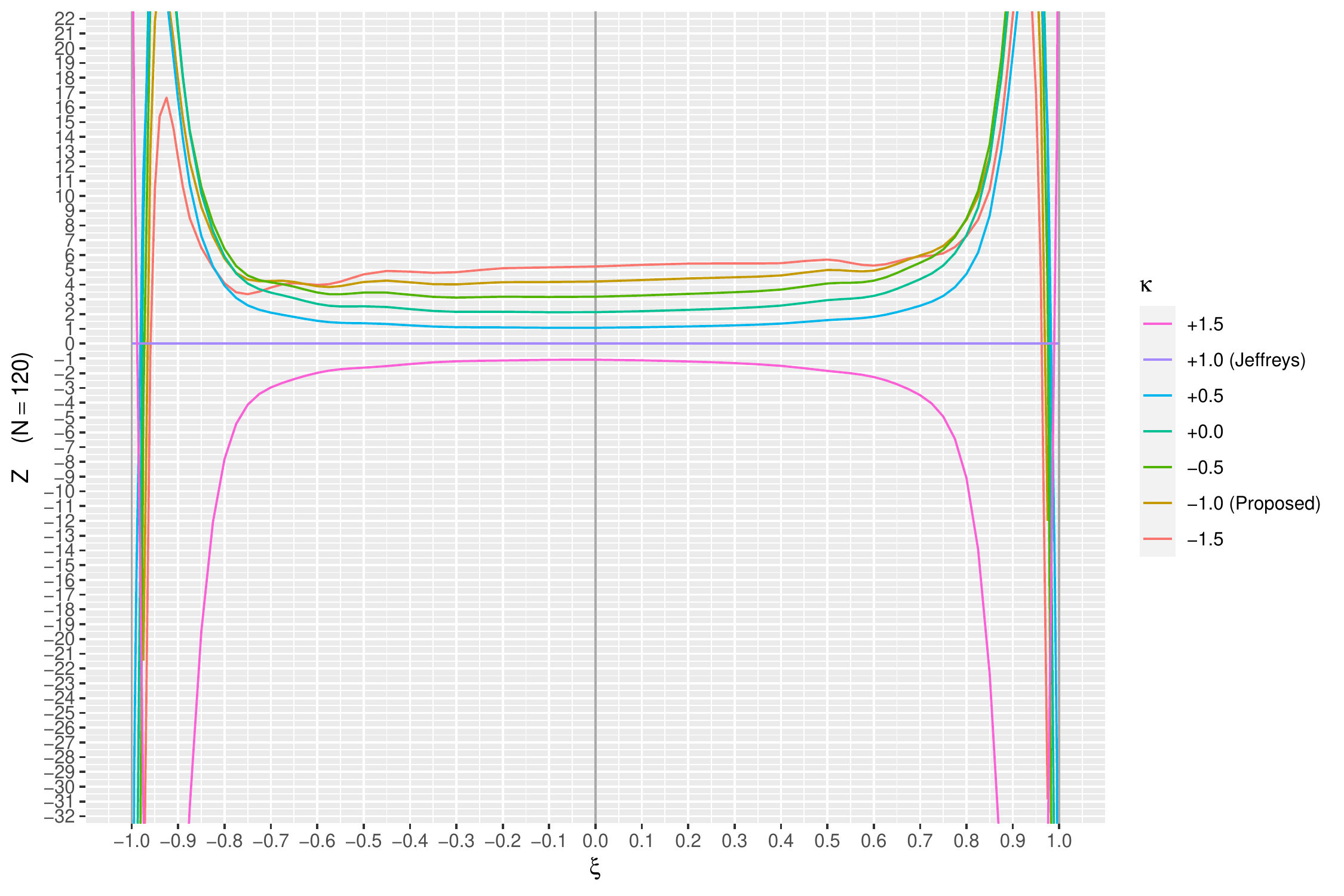}
  \caption{Numerical experiment for the risk difference $Z^{(\kappa)} (\xi)$ for $\mathrm{AR} (1; \mathbb{C})$ with $N = 120$, where $\xi \in (-1, 1) \subset U$. The proposed $\hat{S}^{(N)}_{\pi^{(-1)}}$ dominates the baseline $\hat{S}^{(N)}_{\pi_J}$ for $\xi \in [-0.960, +0.960]$.}
  \label{figure:KL_DIFF_N120}
\end{figure}

The results of numerical experiments for $Z^{(\kappa)} (\xi)$ with $N = 30$, $N = 60$, and $N = 120$ are shown in Figure \ref{figure:KL_DIFF_N30}, \ref{figure:KL_DIFF_N60}, and \ref{figure:KL_DIFF_N120}, respectively,
where the Monte Carlo method is used for evaluating the value of (\ref{definition:BPD}).
These Figures show that $Z^{(\kappa)} (\xi)$ asymptotically achieves $Q^{(\kappa)} (\xi)$, but the rate of convergence may depend on $\xi \in U$.
It appears that the convergence is not uniform on $U$.

\clearpage

\appendices

For simplicity, we intentionally omitted mathematical details for some notations in the preceding sections.
In Appendix \ref{section:wirtinger_calculus}, we introduce complex differential analysis and define the Laplacian as the differential operator on a K\"ahler parameter space.
In Appendix \ref{section:kl_divergence}, we derive and justify the risk for a predictive power spectral density for a complex-valued Gaussian process.
In Appendix \ref{section:hermite_polynomial}, we present the complexified version of tensorial Hermite polynomials.
In Appendix \ref{proof:expectation_loglikehood}, we relate the theory of power spectral densities to the theory of probability densities by utilizing the asymptotic expansions of the expectation of the derivatives of the log likelihood.
In Appendix \ref{section:AEP_EPD}, we present the asymptotic expansion of the risk of the estimative power spectral density with the maximum likelihood estimator.
In Appendix \ref{section:AEP_BPD}, we provide the asymptotic expansion of the Bayesian predictive power spectral density around the maximum likelihood estimator.
In Appendix \ref{section:existenece_BPSD}, we prove the existence of Bayesian predictive power spectral densities for $\mathrm{AR} (p; \mathbb{C})$.
In Appendix \ref{proof:main}, we prove the Main Theorem.
In Appendix \ref{section:alpha_parallel}, we show the relation between $\kappa$-priors and $\alpha$-parallel priors.

\section{Wirtinger calculus}
\label{section:wirtinger_calculus}

In this appendix, we introduce an elegant equivalent formulation of usual differential calculus, called Wirtinger calculus or $\mathbb{C}\mathbb{R}$ calculus \cite{mandic2009complex}.
Let us consider a complex-valued function $f$ defined on $\mathbb{C}^p$.
A function defined on the domain $\mathbb{C}^p$ is always regarded as one defined on the domain $\mathbb{R}^{2p}$.
For the $i$-th complex coordinate $z^i = x^i + \sqrt{-1} \, y^i$ in $\mathbb{C}^p$, the Wirtinger derivatives $\pHd{i}$ and $\pAd{i}$ are defined as linear partial differential operators of the first order:
\begin{align}
  \pHd{i} := \frac{1}{2} \left( \frac{\partial}{\partial x^{i}} - \sqrt{-1} \frac{\partial}{\partial y^{i}} \right)
  \;,\;
  \pAd{i} := \frac{1}{2} \left( \frac{\partial}{\partial x^{i}} + \sqrt{-1} \frac{\partial}{\partial y^{i}} \right)
  \,,
  \label{definition:wirtinger_derivative} 
\end{align}
where $\frac{\partial}{\partial x^{i}}$ and $\frac{\partial}{\partial y^{i}}$ denote the usual partial differential operators on $\mathbb{R}^{2p}$.
The symbol $\pAd{i}$ is occasionally expressed as $\bar{\partial}_i$.

We should mention that although the variables $z^i$ and $\overline{z^i}$ are not independent,
the derivatives $\pHd{i}$ and $\pAd{i}$ are independent differential operators in the complexified tangent space of $\mathbb{C}^p = \mathbb{R}^{2p}$.
In fact, direct computation shows that the set $\{ \pHd{1}, \cdots, \pHd{p}, \pAd{1}, \cdots, \pAd{p} \}$ forms a basis of the complexified tangent space of $\mathbb{C}^p = \mathbb{R}^{2p}$.

The Wirtinger derivatives are not the partial derivatives in usual differential calculus;
however, Wirtinger calculus inherits most of the properties of usual differential calculus.

The most fascinating property inherited by Wirtinger calculus from usual differential calculus is its chain rule property:
\begin{align}
  \pHd{i} (f \circ g) &= \sum_{j = 1}^q (\pHd{i} g^j) (\pHd{j} f \circ g ) + \sum_{j = 1}^q (\pHd{i} \bar{g}^{j}) (\pAd{j} f \circ g ) 
  ,
  \label{formula:chain_rule_h} \\
  \pAd{i} (f \circ g) &= \sum_{j = 1}^q (\pAd{i} g^j) (\pHd{j} f \circ g ) + \sum_{j = 1}^q (\pAd{i} \bar{g}^{j}) (\pAd{j} f \circ g )
  \label{formula:chain_rule_a} 
\end{align}
for $i = 1, \cdots, p$,
where $f \colon \mathbb{C}^q \to \mathbb{C}$ and $g = (g^1, \cdots, g^q) \colon \mathbb{C}^p \to \mathbb{C}^q$.

Another important property of Wirtinger calculus is its summation rule.
For a complex vector $\lambda = (\lambda^1, \cdots, \lambda^p) \in \mathbb{C}^p$ and a complex-valued function $f$ defined on $\mathbb{C}^p$,
\begin{align}
  \sum_{i = 1}^p \lambda^i \pHd{i} f + \sum_{j = 1}^{p} \bar{\lambda}^{j} \pAd{j} f
  = \sum_{k = 1}^p \Re(\lambda^k) \frac{\partial f}{\partial x^k} + \sum_{k = 1}^p \Im(\lambda^k) \frac{\partial f}{\partial y^k}
  \label{forumula:summation_rule},
\end{align}
where the function $f$ is regarded as a function defined on $\mathbb{R}^{2p}$ on the right-hand side of the equation.

As the Einstein notation is assumed throughout this paper,
the summation is automatically taken over those indices that appear exactly twice, once as a superscript and once as a subscript.
Therefore, when the Einstein notation is used, the left-hand side of (\ref{forumula:summation_rule}) is denoted by $\lambda^{\alpha} \pd{\alpha} f$, if $\alpha$ runs through the indices $\{1, \cdots, p, \bar{1}, \cdots, \bar{p} \}$, or sometimes by $\lambda^{i} \pHd{i} f + \lambda^{\bar{j}} \pAd{j} f$, if $i, j$ run through the indices $\{ 1, \cdots, p \}$, where $\lambda^{\bar{j}}$ represents the complex conjugate of $\lambda^j$.
In this paper, we attempt to use the symbols $\alpha, \beta, \gamma, \cdots$ when they run through the indices $\{1, \cdots, p, \bar{1}, \cdots, \bar{p} \}$ and to use the symbols $i, j, k, \cdots$ when they run through the indices $\{1, \cdots, p\}$.

Consider a positive definite metric $\gdd{\alpha}{\beta}$ on $\mathbb{C}^p$, i.e.,
$\gdd{\alpha}{\beta} = \gdd{\beta}{\alpha}$
for any $\alpha, \beta \in \{ 1, \cdots, p, \bar{1}, \cdots, \bar{p} \}$,
and
\begin{align}
  \gdd{\alpha}{\beta} \lambda^{\alpha} \lambda^{\beta}
  = \gHdHd{i}{j} \lambda^{i} \lambda^{j}
  + \gHdAd{i}{j} \lambda^{i} \lambda^{\bar{j}}
  + \gAdHd{i}{j} \lambda^{\bar{i}} \lambda^{j}
  + \gAdAd{i}{j} \lambda^{\bar{i}} \lambda^{\bar{j}}
  > 0
\end{align}
for any $\lambda \in \mathbb{C}^p \setminus \{0\}$.
Let us denote by $\begin{bmatrix} \guu{\alpha}{\beta} \end{bmatrix}$ the inverse matrix of the $2p \times 2p$ matrix $\begin{bmatrix} \gdd{\alpha}{\beta} \end{bmatrix}$,
i.e., $\gdd{\alpha}{\gamma} \guu{\gamma}{\beta}  = {\delta_{\alpha}}^{\beta}$ for the Kronecker delta ${\delta_{\alpha}}^{\beta}$.

For the $i$-th complex coordinate $z^i = x^i + \sqrt{-1} \, y^i$ in $\mathbb{C}^p$, the derivatives $\pHu{i}$ and $\pAu{i}$ are defined by
\begin{align}
  \pHu{i} &:= \gHuHu{i}{j} \pHd{j} + \gHuAu{i}{j} \pAd{j}
  \;,\; \\
  \pAu{i} &:= \gAuHu{i}{j} \pHd{j} + \gAuAu{i}{j} \pAd{j} 
  \;.
\end{align}
The derivatives $\pHu{i}$ and $\pAu{i}$ are also simply defined by $\pu{\alpha} = \guu{\alpha}{\beta} \pd{\beta}$ for $\alpha \in \{1, \cdots, p, \bar{1}, \cdots, \bar{p} \}$.
The symbol $\pAu{i}$ is sometimes denoted by $\bar{\partial}^i$.
For later use, let us define the differential operators $D_{\alpha_1 \cdots \alpha_a}$ and $D^{\alpha_1 \cdots \alpha_a}$ by
\begin{align}
  &D_{\alpha_1 \cdots \alpha_a} := \pd{\alpha_1} \cdots \pd{\alpha_a} \,, \\
  &D^{\alpha_1 \cdots \alpha_a} := \guu{\alpha_1}{\beta_1} \cdots \guu{\alpha_a}{\beta_a} \, \pd{\beta_1} \cdots \pd{\beta_a} \,,
  \label{definition:D_a}
\end{align}
respectively, for $\alpha_1, \cdots, \alpha_a \in \{ 1, \cdots, p, \bar{1}, \cdots, \bar{p} \}$.
In particular, $D_{\alpha} = \pd{\alpha}$ and $D^{\alpha} = \pu{\alpha}$.

A metric $\gdd{\alpha}{\beta}$ is called Hermitian, if
\begin{align}
  \gHdHd{i}{j} = \gAdAd{i}{j} = 0 
  \;,\;
  \gHdAd{i}{j} = \gAdHd{j}{i} = \overline{\gHdAd{j}{i}} = \overline{\gAdHd{i}{j}}
  \label{definition:hermitian_manifold}
\end{align}
for all $i, j = 1, \cdots, p$; see Section 8.4 in \cite{nakahara2003geometry}.
If the metric $\gdd{\alpha}{\beta}$ is Hermitian, the square of the distance of the infinitesimal complex vector $ds$ is given by
\begin{align}
  ds^2 = \gdd{\alpha}{\beta} \, ds^{\alpha} ds^{\beta} = 2 \, \gHdAd{i}{j} \, ds^{i} ds^{\bar{j}}
  .
\end{align}
If the metric $\gdd{\alpha}{\beta}$ is Hermitian, we only need to consider one-fourth of the $2p \times 2p$ complex-valued matrix $\begin{bmatrix} \gdd{\alpha}{\beta} \end{bmatrix}$, namely the $p \times p$ Hermitian matrix $\begin{bmatrix} \gHdAd{i}{j} \end{bmatrix}$.
A complex manifold with a Hermitian metric is called a Hermitian manifold.
The Jacobi formula for the Hermitian manifold is
\begin{align}
  \label{formula:jacobi}
  \gHuAu{i}{j} \pHd{k} \gHdAd{i}{j} = \pHd{k} \log \pi_J
  \;,\;
  \gHuAu{i}{j} \pAd{k} \gHdAd{i}{j} = \pAd{k} \log \pi_J,
\end{align}
where $\pi_J$ is the determinant of the $p \times p$ Hermitian matrix $\begin{bmatrix} \gHdAd{i}{j} \end{bmatrix}$, i.e., the square root of the determinant of the $2p \times 2p$ matrix $\begin{bmatrix} \gdd{\alpha}{\beta} \end{bmatrix}$; see Section 8.4 in \cite{nakahara2003geometry}.

A Hermitian manifold with a metric $\gHdAd{i}{j}$ is called a K\"ahler manifold, if
\begin{align}
  \pHd{i} \gHdAd{j}{k} = \pHd{j} \gHdAd{i}{k}
  \;,\;
  \pAd{i} \gHdAd{j}{k} = \pAd{k} \gHdAd{j}{i}
  \label{definition:kaehler_manifold}
\end{align}
for all $i, j, k = 1, \cdots, p$; see Section 8.5 in \cite{nakahara2003geometry}.

The Laplacian (Laplace--Beltrami operator) on a K\"ahler manifold is
\begin{align}
  \Delta = \pu{\alpha} \pd{\alpha} = \guu{\alpha}{\beta} \pd{\alpha} \pd{\beta} = 2 \gHuAu{i}{j} \pHd{i} \pAd{j}
  \label{formula:laplacian} ,
\end{align}
which does not hold in general for the usual Riemannian manifold.
From (\ref{formula:laplacian}), we have
\begin{align}
  \frac{\Delta \phi^a}{\phi^a} = a \frac{\Delta \phi}{\phi} + 2 a (a - 1) \gHuAu{i}{j} (\pHd{i} \log \phi)(\pAd{j} \log \phi)
  \label{formula:laplacian_kappa}
\end{align}
for $a \in \mathbb{R}$,
which is a useful formula for calculating $\Delta \phi^a$.

\section{Kullback--Leibler divergence between power spectral densities}
\label{section:kl_divergence}

This appendix derives and justifies the form (\ref{definition:risk}) of risk $R ( \hat{S} \mid \theta )$ for a predictive power spectral density $\hat{S}$ for a complex-valued Gaussian process, and it explains why the Bayesian predictive power spectral density $\hat{S}^{(N)}_{\pi}$ minimizes the Bayes risk (\ref{definition:bayes_risk}) given an observation $z^{(N)} \in \mathbb{C}^N$ and a prior $\pi$.

Let $W^{(N)}$ be a circulant matrix and $D^{(N)}$ be a diagonal matrix defined by
\begin{align*}
  &W^{(N)} := \begin{pmatrix}
    0 & 1 & \cdots & 0 & 0 \\
    0 & 0 & \cdots & 0  & 0 \\
    \vdots & \vdots &  \ddots & \vdots & \vdots \\
    0 & 0  & \cdots  & 0  & 1 \\
    1 & 0  &  \cdots & 0 & 0 
  \end{pmatrix}
  ,\\
  &D^{(N)} := \begin{pmatrix}
    e^{2 \pi i \frac{1}{N}} & 0 & \cdots & 0 \\
    0 & e^{2 \pi i \frac{2}{N}} & \cdots & 0 \\
    \vdots & \vdots & \ddots & \vdots \\
    0 & 0 & \cdots & e^{2 \pi i \frac{N}{N}}
  \end{pmatrix}
  ,
\end{align*}
respectively.
We have a relation $W^{(N)} U^{(N)} = U^{(N)} D^{(N)}$,
where $U^{(N)} := \begin{bmatrix} U^{(N)}_{st} \end{bmatrix}$ is a unitary matrix defined by $U^{(N)}_{st} := \frac{1}{\sqrt{N}} e^{2 \pi i \frac{st}{N}}$.
Set $\Lambda (z) := \sum_{h = - \infty}^{\infty} \gamma_h z^{-h}$, where $z$ is a formal variable.
Note that $S(\omega) = \frac{1}{2 \pi} \Lambda(e^{i\omega})$ and $\Lambda \bigl( W^{(N)} \bigr) = U^{(N)} \, \Lambda \bigl( D^{(N)} \bigr) \, { U^{(N)} }^*$.
Suppose $\Lambda (z)$ has a Laurent expansion on $|z| = 1$.
Then, $\bigl( \Lambda (z) \bigr)^{-1} = 1 / \Lambda(z)$ is defined on the neighborhood of $|z| = 1$.
If the autocovariance $\gamma_t$ decreases exponentially, we have $\Lambda \bigl( W^{(N)} \bigr) \approx \Sigma^{(N)}$ for large $N$
because the $(s, t)$-th element of the matrix $\Lambda \bigl( W^{(N)} \bigr)$ is approximated as
\begin{align}
    \begin{bmatrix} \Lambda \bigl( W^{(N)} \bigr) \end{bmatrix}_{st} = \sum_{k=-\infty}^{\infty} \gamma_{(s-t) + kN} \;\approx\; \gamma_{(s-t)} = \begin{bmatrix} \Sigma^{(N)} \end{bmatrix}_{st}
    .
\end{align}
With this approximation, the log likelihood
\begin{align}
  l^{(N)} (z^{(N)})
  = -N \log \pi - \log \bigl| \Sigma^{(N)} \bigr| - z^* \bigl( \Sigma^{(N)} \bigr)^{-1} z
  \label{definition:log_likelihood}
\end{align}
of the observation $z = z^{(N)}$ from a complex-valued Gaussian process with mean $0 \in \mathbb{C}$ is approximated as
\begin{align}
  l^{(N)} (z^{(N)})
  &\approx -N \log \pi - \log | \Lambda(W) | - z^* \bigl( \Lambda (W) \bigr)^{-1} z
  \nonumber\\
  &= N C - \sum_{n=1}^N \log S \left(2 \pi \frac{n}{N} \right) - \sum_{n=1}^N \frac{I \left(2 \pi \frac{n}{N} \right)}{S \left(2 \pi \frac{n}{N} \right)}
  ,\label{formula:approximation_log_likelihood}
\end{align}
where $I$ denotes the empirical power spectral density (periodogram) defined by
$ I \left(2 \pi \frac{n}{N} \right) := \frac{1}{2\pi} | \tilde{z}_n |^2$ 
with $\tilde{z}_n := (U^* z)_n = \frac{1}{\sqrt{N}} \sum_{s=1}^N e^{- 2 \pi i s \frac{n}{N}} z_s$ and $C$ is a constant independent of $N$ and $S$.
For a more detailed explanation of (\ref{formula:approximation_log_likelihood}),
see \cite{whittle1953estimation} for real-valued stationary processes and \cite{anderson1977estimation} for real-valued ARMA processes.

Suppose the variance-covariance matrix (\ref{definition:variance_covariance_matrix}) is parameterized by complex parameters $\theta \in \Theta \subset \mathbb{C}^p$,
i.e., the autocovariances $\{ \gamma_h \}_{h \in \mathbb{Z}}$ are parametrized by $\theta \in \Theta$.
Its power spectral density (\ref{definition:power_spectral_density}) is denoted by $S_{\theta} (\omega)$ or $S (\omega \mid \theta)$ for $\theta \in \Theta$.
For $\theta \in \Theta$, we denote the corresponding probability distribution, probability density function, and log likelihood of the observation $z^{(N)} \in \mathbb{C}^N$ by $P^{(N)}_{\theta}$, $p^{(N)}_{\theta}$, and $l^{(N)}_{\theta}$, respectively.

For $\theta_1, \theta_2 \in \Theta$, the KL-divergence $\KL{P_{\theta_1}}{P_{\theta_2}}$ of the distributions $P_{\theta_2}$ from the distribution $P_{\theta_1}$ is approximated as
\begin{align*}
  \KL{P_{\theta_1}}{P_{\theta_2}}
  &= \int_{\mathbb{C}^p} \left(
    l^{(N)}_{\theta_1} \bigl( z^{(N)} \bigr)
    - l^{(N)}_{\theta_2} \bigl( z^{(N)} \bigr)
  \right) dP_{\theta_1}  \bigl( z^{(N)} \bigr) \\
  &\approx E_{\theta_1}
  \left[
    - \sum_{n=1}^N \log \frac{S_{\theta_1} \left(2 \pi \frac{n}{N} \right)}{S_{\theta_2} \left(2 \pi \frac{n}{N} \right)}
    - \sum_{n=1}^N \left( \frac{I \left(2 \pi \frac{n}{N} \right)}{S_{\theta_1} \left(2 \pi \frac{n}{N} \right)} - \frac{I \left(2 \pi \frac{n}{N} \right)}{S_{\theta_2} \left(2 \pi \frac{n}{N} \right)} \right)
  \right] \\
  &\approx \frac{N}{2\pi} \int_{-\pi}^{\pi} \left( - \log \frac{S_{\theta_1} (\omega)}{S_{\theta_2} (\omega)} - 1 + \frac{S_{\theta_1} (\omega)}{S_{\theta_2} (\omega)} \right) d\omega \\
  &= N \KL{S_{\theta_1}}{S_{\theta_2}}
  ,
\end{align*}
where $\KL{S_{\theta_1}}{S_{\theta_2}}$ is the KL-divergence (\ref{definition:kl_divergence}) between power spectral densities, which has previously been discussed in the literature \cite{amari1987differential, choi2015kahlerian} in the context of signal processing.
On the other hand, in the literature \cite{komaki1999estimating, tanaka2011asymptotic} on real-valued processes, $4 \pi$ instead of $2 \pi$ is used in the denominator of (\ref{definition:kl_divergence}).
However, as we have explained, the constant in the denominator of (\ref{definition:kl_divergence}) for the complex-valued process should be $2 \pi$.
Note also that for any power spectral densities $S_1$ and $S_2$, because $- \log x -1 + x \geq 0$ for any $x \geq 0$, we have $\KL{S_1}{S_2} \geq 0$ in general, and $\KL{S_1}{S_2} = 0$ if and only if $S_1 (\omega) = S_2 (\omega)$ for $\omega \in [-\pi, \pi]$ almost everywhere.
The asymptotic expansion
\begin{align}
  - \log \frac{1}{1 + x} - 1 + \frac{1}{1 + x}
  = \frac{1}{2} x^2 - \frac{2}{3} x^3 + \frac{3}{4} x^4 - \frac{4}{5} x^5 + O (x^6)
  \label{formula:asymptotic_expansion_kl_divergence}
\end{align}
is a useful formula for calculating the value of $\KL{S}{S + dS}$.

For a possibly improper prior $\pi$, the Bayesian predictive power spectral density $\hat{S}^{(N)}_{\pi}$ minimizes the Bayes risk (\ref{definition:bayes_risk}) among all the predictive power spectral densities $\hat{S}^{(N)}$ if $r \bigl( \hat{S}_{\pi}^{(N)} \mid \pi \bigr) < + \infty$; see \cite{aitchison1975goodness}.
In fact,
\begin{align*}
  r \bigl( \hat{S}^{(N)} \mid \pi \bigr) - r \bigl( \hat{S}^{(N)}_{\pi} \mid \pi \bigr)
  &= \int_{\Theta} \int_{\mathbb{C}^p} \left(
    \KL{ S_{\theta} }{ \hat{S}^{(N)}} - \KL{ S_{\theta} }{ \hat{S}_{\pi}^{(N)}}
  \right) dP_{\theta}^{(N)} d\Theta \\
  &= \int_{\Theta} \int_{\mathbb{C}^p}
  \left(
    \frac{1}{2 \pi} \int_{-\pi}^{\pi} \left(
      - \log \frac{ \hat{S}^{(N)}_{\pi} }{ \hat{S}^{(N)} } - \frac{S_{\theta}}{\hat{S}^{(N)}_{\pi}} + \frac{S_{\theta}}{\hat{S}^{(N)}}
    \right) d\omega
  \right) dP_{\theta}^{(N)} d\Theta
  \\
  &= \int_{\mathbb{C}^p} \KL{\hat{S}^{(N)}_{\pi}}{\hat{S}^{(N)}} m^{(N)}_{\pi} \bigl( z^{(N)} \bigr) \,dz^{(N)}
  \geq 0
\end{align*}
for any predictive power spectral density $\hat{S}^{(N)}$,
where $dP_{\theta}^{(N)} := p_{\theta}^{(N)} \bigl( z^{(N)} \bigr) dz^{(N)}$, $d\Theta := \pi(\theta) d\theta$,
and
$ m^{(N)}_{\pi} \bigl( z^{(N)} \bigr) := \int_{\Theta} p^{(N)}_{\theta} \bigl( z^{(N)} \bigr) \, \pi(\theta) \, d\theta $
is the marginal distribution of the observation $z^{(N)} \in \mathbb{C}^N$ based on the prior $\pi$.

\section{Tensorial Hermite polynomials}
\label{section:hermite_polynomial}

Tensorial Hermite polynomials, as introduced in \cite{amari1983differential}, are very useful tools for calculating Edgeworth expansions.
Here, we present the complexified version of tensorial Hermite polynomials to calculate (\ref{formula:AE_I}).

Consider a metric $\gdd{\alpha}{\beta}$ on $\mathbb{C}^p$.
We define a complex-valued function $\phi$ on $\mathbb{C}^p$ by
$ \phi(\lambda) := \frac{1}{G} \, e^{- \frac{1}{2} \gdd{\alpha}{\beta} \lambda^{\alpha} \lambda^{\beta}} $
for $\lambda = (\lambda^1, \cdots, \lambda^p) \in \mathbb{C}^p$,
where $\alpha, \beta$ run through the indices $\{ 1, \cdots, p, \bar{1}, \cdots, \bar{p} \}$
and
$ G := \int_{\mathbb{C}^p} e^{- \frac{1}{2} \gdd{\alpha}{\beta} \lambda^{\alpha} \lambda^{\beta}} d\lambda $
is the normalization factor that gives $\int_{\mathbb{C}^p} \phi(\lambda) \,d\lambda = 1$.
We assume that the metric $\gdd{\alpha}{\beta}$ is positive definite so that $\lim_{|\lambda| \to \infty} |\phi(\lambda)| = 0$.
If the metric $\gdd{\alpha}{\beta}$ is Hermitian,
the normalization factor $G$ reduces to the product of $\pi^p$ and the determinant of the $p \times p$ Hermitian matrix
$\begin{bmatrix} \gHuAu{i}{j} \end{bmatrix}$.
However, to obtain a general result, we do not assume that the metric $\gdd{\alpha}{\beta}$ is Hermitian in this appendix.

We define the complex-valued tensorial Hermite polynomial $h^{\alpha_1 \cdots \alpha_a}$ for $\alpha_1, \cdots, \alpha_a \in \{ 1, \cdots, p, \bar{1}, \cdots, \bar{p} \}$ by the identity
\begin{align}
  (-1)^k D^{\alpha_1 \cdots \alpha_a} \phi (\lambda) = h^{\alpha_1 \cdots \alpha_a} (\lambda) \, \phi (\lambda)
  ,
\end{align}
where the differential operator $D^{\alpha_1 \cdots \alpha_a}$ is defined by (\ref{definition:D_a}).
For example, $h^{\alpha_1} (\lambda) = \lambda^{\alpha_1}$, and $h^{\alpha_1 \alpha_2} (\lambda) = \lambda^{\alpha_1} \lambda^{\alpha_2} - \guu{\alpha_1}{\alpha_2}$.
Following a procedure similar to that in in \cite{amari1983differential},
we obtain
\begin{align*}
  &\int_{\mathbb{C}^p} h^{\alpha_1 \cdots \alpha_a} (\lambda) \, h^{\beta_1 \cdots \beta_b} (\lambda) \, \phi (\lambda) \,d\lambda
  =
  \begin{cases}
    a! \, \guu{(\alpha_1}{\beta_1} \cdots \guu{\alpha_a}{\beta_a)} &\quad\text{if $a = b$} \\
    0 &\quad\text{otherwise}
  \end{cases}
\end{align*}
for $\alpha_1, \cdots, \alpha_a, \beta_1, \cdots, \beta_b \in \{ 1, \cdots, p, \bar{1}, \cdots, \bar{p} \}$,
where the parentheses around the indices $\alpha_1 \beta_1 \cdots \alpha_a \beta_a$ imply the symmetrization of the indices $\alpha_1, \cdots, \alpha_a$;
i.e., $(\alpha_1 \beta_1 \cdots \alpha_a \beta_a) := (1 / a!) \sum_{\sigma} \alpha_{\sigma(1)} \beta_1 \cdots \alpha_{\sigma(a)} \beta_a$, where $\sigma$ runs through all the permutations of the indices $\alpha_1, \cdots, \alpha_a$.
For example, 
\begin{align}
  \int_{\mathbb{C}^p} \lambda^{\alpha_1} \lambda^{\alpha_2} \,\phi(\lambda) \,d\lambda
  = \int_{\mathbb{C}^p} (h^{\alpha_1 \alpha_2} + \guu{\alpha_1}{\alpha_2}) \,\phi(\lambda) \,d\lambda = \guu{\alpha_1}{\alpha_2},
  \label{formula:AE_I2}
\end{align}
and
\begin{align}
  \int_{\mathbb{C}^p} \lambda^{\alpha_1} \lambda^{\alpha_2} \lambda^{\alpha_3} \lambda^{\alpha_4} \,\phi(\lambda) \,d\lambda
  &= \int_{\mathbb{C}^p} (
    h^{\alpha_1 \alpha_2} h^{\beta_1 \beta_2}
    + \lambda^{\alpha_1} \lambda^{\alpha_2} \guu{\beta_1}{\beta_2}
    + \guu{\alpha_1}{\alpha_2} \lambda^{\beta_1} \lambda^{\beta_2}
    - \guu{\alpha_1}{\alpha_2} \guu{\beta_1}{\beta_2}
  ) \,\phi(\lambda) \,d\lambda \nonumber\\
  &= \guu{\alpha_1}{\beta_1} \guu{\alpha_2}{\beta_2} + \guu{\alpha_2}{\beta_1} \guu{\alpha_1}{\beta_2} + \guu{\alpha_1}{\alpha_2} \guu{\beta_1}{\beta_2}
  \label{formula:AE_I4} 
\end{align}
because $\int h^{\alpha_1 \alpha_2} (\lambda) \, h^{\beta_1 \beta_2} (\lambda) \,\phi(\lambda) \,d\lambda = 2 \, \guu{(\alpha_1}{\beta_1} \guu{\alpha_2}{\beta_2)} = \guu{\alpha_1}{\beta_1} \guu{\alpha_2}{\beta_2} + \guu{\alpha_2}{\beta_1} \guu{\alpha_1}{\beta_2}$.

\section{Asymptotic expansion of the expectation of the derivatives of the log likelihood}
\label{proof:expectation_loglikehood}

Here, we show the asymptotic expansions of the expectation of the derivatives of the log likelihood.
These asymptotic expansions relate the theory of power spectral densities to the theory of probability densities.

Consider a variance-covariance matrix $\Sigma$ of size $N$ parametrized by complex parameters $\theta = (\theta^1, \cdots, \theta^p) \in \Theta \subset \mathbb{C}^p$.
We set the matrix
\begin{align*}
  \Sigma_{\alpha_1 \cdots \alpha_a , \beta_1 \cdots \beta_b, \cdots, \gamma_1 \cdots \gamma_c} 
  &:=
  \Sigma^{-1} \left( D_{\alpha_1 \cdots \alpha_a} \Sigma \right)
  \Sigma^{-1} \left( D_{\beta_1 \cdots \beta_b} \Sigma \right)
  \cdots \Sigma^{-1} \left( D_{\gamma_1 \cdots \gamma_c} \Sigma \right)
  .
\end{align*}
For example, $\Sigma_{1\bar{1},2} = \Sigma^{-1} (\pHd{1} \pAd{1} \Sigma) \Sigma^{-1} (\pHd{2} \Sigma)$.
Direct computation shows that
\begin{align*}
  \pd{\alpha} l =&\,
  z^* ( \Sigma_{\alpha} ) \Sigma^{-1} z - \tr \left( \Sigma_{\alpha} \right)
  ,\\
  \pd{\alpha} \pd{\beta} l =& 
  - z^* \left(
    \Sigma_{\alpha, \beta}
    + \Sigma_{\beta, \alpha}
  \right) \Sigma^{-1} z
  + z^* \left(
    \Sigma_{\alpha \beta}
  \right) \Sigma^{-1} z
  - \tr
  \left(
    - \Sigma_{\alpha, \beta}
    + \Sigma_{\alpha \beta}
  \right)
  ,\\
  \pd{\alpha} \pd{\beta} \pd{\gamma} l
  =&\, z^* \left(
    \Sigma_{\alpha, \beta, \gamma}
    + \Sigma_{\alpha, \gamma, \beta}
    + \Sigma_{\beta, \alpha, \gamma}
    + \Sigma_{\beta, \gamma, \alpha}
    + \Sigma_{\gamma, \alpha, \beta}
    + \Sigma_{\gamma, \beta, \alpha}
  \right) \Sigma^{-1} z \\
  &- z^* \left(
    \Sigma_{\alpha, \beta \gamma}
    + \Sigma_{\beta, \gamma \alpha}
    + \Sigma_{\gamma, \alpha \beta}
    + \Sigma_{\alpha \beta, \gamma}
    + \Sigma_{\beta \gamma, \alpha}
    + \Sigma_{\gamma \alpha, \beta}
  \right) \Sigma^{-1} z \\
  &+ z^* \left(
    \Sigma_{\alpha \beta \gamma}
  \right) \Sigma^{-1} z 
  -\tr
  \left(
    \Sigma_{\alpha, \beta, \gamma}
    + \Sigma_{\alpha, \gamma, \beta}
    - \Sigma_{\alpha, \beta \gamma}
    - \Sigma_{\beta, \alpha \gamma}
    - \Sigma_{\gamma, \alpha \beta}
    + \Sigma_{\alpha \beta \gamma}
  \right)
  ,
\end{align*}
where $l = l^{(N)}_{\theta} \left( z^{(N)} \right)$ is the log likelihood (\ref{definition:log_likelihood}) of the observation $z^{(N)} \in \mathbb{C}^N$ from the complex-valued Gaussian process at the parameter $\theta \in \Theta  \subset \mathbb{C}^p$.

As $E [ z z^* ] = \Sigma$, we have
\begin{align*}
  \frac{1}{N} \, E_{\theta} \left[ \pd{\alpha} l^{(N)}_{\theta} \right]
  &= 0 \,,\\
  \frac{1}{N} \, E_{\theta} \left[ \pd{\alpha} \pd{\beta} l^{(N)}_{\theta} \right]
  &= - \tr \, \Sigma_{\alpha, \beta} \,, \\
  \frac{1}{N} \, E_{\theta} \left[ \pd{\alpha} \pd{\beta} \pd{\gamma} l^{(N)}_{\theta} \right]
  &= 2 \, \tr \, \bigl( \Sigma_{\alpha, \beta, \gamma}  + \Sigma_{\gamma, \beta, \alpha} \bigr)
  - \tr \, \Sigma_{\alpha \beta, \gamma}  - \tr \, \Sigma_{\beta \gamma, \alpha}  - \tr \, \Sigma_{\gamma \alpha, \beta} \,,
\end{align*}
where $E_{\theta}$ denotes the expectation over the distribution $P^{(N)}_{\theta}$ of the observation $z^{(N)} \in \mathbb{C}^N$ at the parameter $\theta \in \Theta$. 

Utilizing complex-valued tensorial polynomials, we also have
\begin{align*}
  \frac{1}{N} \, E_{\theta} \left[ \bigl( \pd{\alpha} l^{(N)}_{\theta} \bigr) \bigl( \pd{\beta} l^{(N)}_{\theta} \bigr) \right]
  &= \tr \, \Sigma_{\alpha, \beta} \,, \\
  \frac{1}{N} \, E_{\theta} \left[ \bigl( \pd{\alpha} \pd{\beta} l^{(N)}_{\theta} \bigr) \bigl( \pd{\gamma} l^{(N)}_{\theta} \bigr) \right]
  &= \tr \, \Sigma_{\alpha \beta, \gamma} 
  - \, \tr \, \bigl( \Sigma_{\alpha, \beta, \gamma} + \Sigma_{\gamma, \beta, \alpha} \bigr) \,,
\end{align*}
and
\begin{align*}
  &\frac{1}{N} \, E_{\theta} \left[ \bigl( \pd{\alpha} l^{(N)}_{\theta} \bigr) \bigl( \pd{\beta} l^{(N)}_{\theta} \bigr) \bigl( \pd{\gamma} l^{(N)}_{\theta} \bigr) \right] 
  = \tr \, \bigl( \Sigma_{\alpha, \beta, \gamma} + \Sigma_{\gamma, \beta, \alpha} \bigr) \,.
\end{align*}

To relate the expectation of the derivative of the log likelihood to the quantities $M_{\alpha_1 \cdots \alpha_a,\, \beta_1 \cdots \beta_b,\, \cdots,\, \gamma_1 \cdots \gamma_c}$ defined as (\ref{definition:M}), we utilize the theorem proved in \cite{taniguchi1983second},
which was originally proved for real-valued processes but is still valid for complex-valued processes.
We introduce the space $\mathcal{D}$ of power spectral densities of complex-valued processes defined on $[-\pi, \pi]$:
\begin{align*}
  \mathcal{D} := \bigl\{ S \;\mid\; S(\omega) = \sum_{h = -\infty}^{\infty} \gamma_h \, e^{- \sqrt{-1}\, h \omega} \,,
  \overline{\gamma_h} = \gamma_{-h} \,,\, \sum_{h = -\infty}^{\infty} |h| |\gamma_h| < \infty \bigr\}
  .
\end{align*}
The space $\mathcal{D}_{\mathrm{ARMA}}$ of power spectral densities for complex-valued stationary ARMA processes, where we have assumed causality and invertibility of the process, is a subspace of $\mathcal{D}$.
\begin{prop}[\cite{taniguchi1983second}] \label{thm:T1983}
  For $S_1, \cdots, S_a \in \mathcal{D}_{\mathrm{ARMA}}$ and $F_1, \cdots, F_a \in \mathcal{D}$,
  \begin{align*}
    \frac{1}{N}\, \tr \left( \Sigma(F_1) \, \Sigma(S_1)^{-1} \cdots \Sigma(F_a) \, \Sigma(S_a)^{-1} \right)
    = \frac{1}{2 \pi} \int_{-\pi}^{\pi} F_1 (\omega) \cdots F_a (\omega)\, S_1^{-1} (\omega) \cdots S_a^{-1}  (\omega) \,d\omega
    + O(N^{-1})
    ,
  \end{align*}
  where the $(s, t)$-th element $\begin{bmatrix} \Sigma (F) \end{bmatrix}_{st}$ of the matrix $\Sigma(F)$ of size $N$ is defined as
  $ \begin{bmatrix} \Sigma (F) \end{bmatrix}_{st} := \int_{-\pi}^{\pi} e^{\sqrt{-1}\, (s - t)}  F(\omega) \,d\omega$
    .
\end{prop}

For a power spectral density $S = S_{\theta} \in \mathcal{D}_{\mathrm{ARMA}}$ of a complex-valued stationary ARMA process parameterized by $\theta \in \Theta  \subset \mathbb{C}^p$
and its variance-covariance matrix $\Sigma = \Sigma_{\theta}^{(N)}$ of size $N$ for $z = z^{(N)} \in \mathbb{C}^N$,
the $(s, t)$-th element $\begin{bmatrix} \pd{\alpha_1} \cdots \pd{\alpha_a} \Sigma \end{bmatrix}_{st}$ of the matrix $\pd{\alpha_1} \cdots \pd{\alpha_a} \Sigma$ is calculated as
\begin{align*}
  \begin{bmatrix}
    \pd{\alpha_1} \cdots \pd{\alpha_a} \Sigma 
  \end{bmatrix}_{st}
  = \int_{-\pi}^{\pi} e^{\sqrt{-1}\, (s - t)} \left( \pd{\alpha_1} \cdots \pd{\alpha_a} S(\omega) \right) d\omega
  .
\end{align*}

Thus, we have the following proposition, which relates the quantities $M_{\alpha_1 \cdots \alpha_a,\, \beta_1 \cdots \beta_b,\, \cdots,\, \gamma_1 \cdots \gamma_c}$ defined as (\ref{definition:M}) to the quantities $K_{\alpha_1 \cdots \alpha_a,\, \beta_1 \cdots \beta_b,\, \cdots,\, \gamma_1 \cdots \gamma_c}$ defined as (\ref{definition:K}).
\begin{prop}
  \label{prop:expectation_loglikehood}
  For a complex-valued stationary ARMA process parameterized by $\theta \in \Theta  \subset \mathbb{C}^p$,
  \begin{align}
    &\frac{1}{N} \, E_{\theta} \left[ \pd{\alpha} l^{(N)}_{\theta} \right]
    = 0
    \label{formula:expectation_loglikehood1}, \\
    &\frac{1}{N} \, E_{\theta} \left[ \pd{\alpha} \pd{\beta} l^{(N)}_{\theta} \right]
    = - \gdd{\alpha}{\beta} + O(N^{-1})
    \label{formula:expectation_loglikehood2}, \\
    &\frac{1}{N} \, E_{\theta} \left[ \pd{\alpha} \pd{\beta} \pd{\gamma} l^{(N)}_{\theta} \right]
    = - \Guddd{(m)}{\alpha}{\beta}{\gamma} - \Guddd{(m)}{\beta}{\gamma}{\alpha} - \Guddd{(m)}{\gamma}{\alpha}{\beta}
    + 2 \Tddd{\alpha}{\beta}{\gamma} + O(N^{-1})
    \label{formula:expectation_loglikehood3}, \\
    &\frac{1}{N} \, E_{\theta} \left[ \bigl( \pd{\alpha} l^{(N)}_{\theta} \bigr) \bigl( \pd{\beta} l^{(N)}_{\theta} \bigr) \right]
    = \gdd{\alpha}{\beta} + O(N^{-1}) ,\\
    &\frac{1}{N} \, E_{\theta} \left[ \bigl( \pd{\alpha} \pd{\beta} l^{(N)}_{\theta} \bigr) \bigl( \pd{\gamma} l^{(N)}_{\theta} \bigr) \right]
    = \Guddd{(m)}{\alpha}{\beta}{\gamma} - \Tddd{\alpha}{\beta}{\gamma} + O(N^{-1}) ,\\
    &\frac{1}{N} \, E_{\theta} \left[ \bigl( \pd{\alpha} l^{(N)}_{\theta} \bigr) \bigl( \pd{\beta} l^{(N)}_{\theta} \bigr) \bigl( \pd{\gamma} l^{(N)}_{\theta} \bigr) \right]
    = \Tddd{\alpha}{\beta}{\gamma} + O(N^{-1})
  \end{align}
  for $\alpha, \beta, \gamma \in \{ 1, \cdots, p, \bar{1}, \cdots, \bar{p} \}$,
  where $l^{(N)}_{\theta}$ denotes the log likelihood (\ref{definition:log_likelihood}) of the observation $z = z^{(N)} \in \mathbb{C}^N$ of size $N$
  and $E_{\theta}$ denotes the expectation over the distribution $P^{(N)}_{\theta}$ of the observation $z^{(N)} \in \mathbb{C}^N$ at the parameter $\theta \in \Theta$. 
  The quantities $\gdd{\alpha}{\beta}$, $\Tddd{\alpha}{\beta}{\gamma}$, and $\Guddd{(m)}{\alpha}{\beta}{\gamma}$ on the right-hand side of the equations are defined as (\ref{definition:M_g}), (\ref{definition:M_T}), and (\ref{definition:M_Gm}), respectively.
\end{prop}

\section{Asymptotic expansion of estimative power spectral densities}
\label{section:AEP_EPD}

Here, we preset the asymptotic expansion of the risk of the estimative power spectral density with the maximum likelihood estimator.
The risk is approximately $\frac{p}{N}$ for most predictive power spectral densities with asymptotically efficient estimators when we use $p$ complex parameters.

Let us fix the true parameter $\theta_0 \in \Theta$ and denote the maximum likelihood estimator by $\hat{\theta}^{(N)}$ for a while;
we set $S_0 := S_{\theta_0}$ and $\hat{S}^{(N)} := S_{\hat{\theta}^{(N)}}$.
For $\lambda := \sqrt{N} \bigl( \hat{\theta}^{(N)} - \theta_0 \bigr) = O_P(1)$,
utilizing (\ref{formula:asymptotic_expansion_kl_divergence}) and the Taylor expansion of $\hat{S}^{(N)}$ around $S_0$,
we obtain 
\begin{align}
  \KL{S_0}{\hat{S}^{(N)}}
  = \frac{1}{2 N} \gdd{\alpha}{\beta} \lambda^{\alpha} \lambda^{\beta}
  + \frac{1}{N \sqrt{N}} \left( \frac{1}{2} \Guddd{(m)}{\alpha}{\beta}{\gamma} - \frac{1}{3} \Tddd{\alpha}{\beta}{\gamma} \right) \lambda^{\alpha} \lambda^{\beta} \lambda^{\gamma} 
  + O_P(N^{-2})
  ,
\end{align}
where the quantities $M_{\alpha_1 \cdots \alpha_a,\, \beta_1 \cdots \beta_b,\, \cdots,\, \gamma_1 \cdots \gamma_c}$ appearing on the right-hand side are all evaluated at $\theta_0$. 
In particular,
the asymptotic expansion of the risk of the maximum likelihood estimator evaluated at the true parameter $\theta_0$ is given by
\begin{align}
  R \bigl( \hat{S}^{(N)} \mid \theta_0 \bigr) = E_{\theta_0} \left[ \KL{S_0}{\hat{S}^{(N)}} \right] = \frac{p}{N} + O(N^{-2})
\end{align}
because $E_{\theta_0} \bigl[ \lambda^{\alpha} \lambda^{\beta} \bigr] = \guu{\alpha}{\beta} + O \bigl( N^{-1} \bigr)$, $E_{\theta_0} \bigl[ \lambda^{\alpha} \lambda^{\beta} \lambda^{\gamma} \bigr] = O (N^{-\frac{1}{2}})$, and $\guu{\alpha}{\beta} \gdd{\alpha}{\beta} = 2p$.

\section{Asymptotic expansion of Bayesian predictive power spectral densities}
\label{section:AEP_BPD}

Here, we present an asymptotic expansion of the Bayesian predictive power spectral density $\hat{S}^{(N)}_{\pi}$ of a complex-valued ARMA process around the maximum likelihood estimator $\hat{\theta}^{(N)}$.
This is the first step to obtaining the asymptotic expansion of the risk differences needed in the proof of the Main Theorem.

We follow the original proof \cite{tanaka2011asymptotic} for the real-valued ARMA process.
However, because we consider complex-valued processes, the definitions of some quantities must be slightly modified.
Basically, the proof for the real-valued ARMA process is applied to the proof for the complex-valued ARMA process because a process parameterized by $p$ complex parameters is essentially a process parameterized by $2p$ real parameters.
However, we must keep in mind that the evenness $S( \omega) = S(- \omega)$ for $\omega \in [-\pi, \pi]$ of power spectral densities is not valid for complex-valued processes.
In what follows, we carefully trace the proof for the real-valued process to ensure that the evenness of power spectral densities is nowhere used in the proof.

For the maximum likelihood estimator $\hat{\theta} = \hat{\theta}^{(N)} \left( z^{(N)} \right) = \theta_0 + O_P(N^{-\frac{1}{2}})$ for the observation $z^{(N)} \in \mathbb{C}^N$ of size $N$ from the complex-valued Gaussian process having the true parameter $\theta_0 \in \Theta \subset \mathbb{C}^p$,
the Bayesian predictive power spectral density is expanded as
\begin{align}
  \hat{S}_{\pi}^{(N)} (\omega)
  = S_{\hat{\theta}} (\omega)
  + \frac{1}{\sqrt{N}} \left( \pd{\alpha} S_{\hat{\theta}} (\omega) \right) E^{\pi} \bigl[ \tilde{\lambda}^{\alpha} \bigr]
  \label{formula:AEP_S1}
  + \frac{1}{2 N} \left( \pd{\alpha} \pd{\beta} S_{\hat{\theta}} (\omega) \right) E^{\pi} \bigl[ \tilde{\lambda}^{\alpha} \tilde{\lambda}^{\beta} \bigr]
  + O_P (N^{-\frac{3}{2}})
\end{align}
around the maximum likelihood estimator $\hat{\theta}$,
where
\begin{align}
  E^{\pi} \left[ \tilde{\lambda}^{\alpha_1} \cdots \tilde{\lambda}^{\alpha_a} \right]
  &:= \int_{\Theta} \tilde{\lambda}^{\alpha_1} \cdots \tilde{\lambda}^{\alpha_a}
  \,\pi \bigl( \theta \mid z^{(N)} \bigr)
  \label{formula:AEP_lambda}
\end{align}
for $\tilde{\lambda} = \sqrt{N} \bigl( \theta - \hat{\theta} \bigr)$ and $\alpha_1, \cdots, \alpha_a \in \{1, \cdots, p, \bar{1}, \cdots, \bar{p} \}$.
To complete the asymptotic expansion of (\ref{formula:AEP_S1}) around the maximum likelihood estimator $\hat{\theta}$,
we require the asymptotic expansions of (\ref{formula:AEP_lambda}) for $a = 1, 2$.

Let us fix, for a while, the observation $z^{(N)} \in \mathbb{C}^N$
and denote the maximum likelihood estimator by $\hat{\theta} = \hat{\theta}^{(N)} (z^{(N)})$.
For any $\theta \in \Theta \subset \mathbb{C}^p$ such that $\tilde{\lambda} = \sqrt{N} \bigl( \theta - \hat{\theta} \bigr) = O(1)$,
the asymptotic expansion of
$ \tilde{l}_{\theta}^{(N)} \left( z^{(N)} \right) := l_{\theta}^{(N)} \left( z^{(N)} \right) + \log \pi(\theta) $
around the maximum likelihood estimator $\hat{\theta}$ is calculated as
\begin{align}
  \tilde{l}^{(N)}_{\theta}
  = \tilde{l}^{(N)}_{\hat{\theta}} - \frac{1}{2} J^{(N)}_{\alpha\beta}  \tilde{\lambda}^{\alpha} \tilde{\lambda}^{\beta}
  + \frac{ A^{(N)} \bigl( \tilde{\lambda} \bigr) }{ \sqrt{N} } + O \left( N^{-1} \right)
  , \label{formula:AEP_l}
\end{align}
where $J^{(N)}_{\alpha\beta} := -\frac{1}{N} \pd{\alpha} \pd{\beta} l_{\hat{\theta}}^{(N)}$ and 
\begin{align}
  A^{(N)} \bigl( \tilde{\lambda} \bigr) &:= \frac{1}{6} \left( \frac{1}{N} \pd{\alpha} \pd{\beta} \pd{\gamma} l_{\hat{\theta}}^{(N)} \right) \tilde{\lambda}^{\alpha} \tilde{\lambda}^{\beta} \tilde{\lambda}^{\gamma}
  + \left( \pd{\alpha} \log \pi (\hat{\theta}) \right) \tilde{\lambda}^{\alpha}
  .
\end{align}

By referring to (\ref{formula:AEP_l}) and the procedure in \cite{philippe2002non}, we can expand (\ref{formula:AEP_lambda}) as
\begin{align}
  E^{\pi} \left[ \tilde{\lambda}^{\alpha_1} \cdots \tilde{\lambda}^{\alpha_a} \right]
  &= \int_{\sqrt{N}(\Theta - \hat{\theta})} \tilde{\lambda}^{\alpha_1} \cdots \tilde{\lambda}^{\alpha_a}
  \frac{1}{G^{(N)}} e^{- \frac{1}{2} J^{(N)}_{\alpha\beta} \tilde{\lambda}^{\alpha} \tilde{\lambda}^{\beta}} 
  \left( 1 + \frac{A^{(N)} \bigl( \tilde{\lambda} \bigr) }{ \sqrt{N} } + O_P (N^{-1}) \right) d\tilde{\lambda}
\end{align}
for $\tilde{\lambda} = \sqrt{N} \bigl( \theta - \hat{\theta} \bigr)$,
where
\begin{align*}
  G^{(N)} := \int_{\sqrt{N}(\Theta - \hat{\theta})} e^{- \frac{1}{2} J^{(N)}_{\alpha\beta} \tilde{\lambda}^{\alpha} \tilde{\lambda}^{\beta}} d\tilde{\lambda}
\end{align*}
is a normalization constant.

The terms in (\ref{formula:AEP_lambda}) for $a = 1, 2$ are expanded as
\begin{align}
  E^{\pi} \left[ \tilde{\lambda}^{\alpha} \right]
  &= \frac{1}{6 \sqrt{N}} L_{\beta\gamma\delta} I^{\alpha\beta\gamma\delta}
  + \frac{1}{\sqrt{N}} \, \left( \pd{\beta} \log \pi \bigl( \hat{\theta} \bigr) \right) \, I^{\alpha\beta}
  + O_P (N^{-1}),
  \label{formula:AEP_lambda1} \\
  E^{\pi} \left[ \tilde{\lambda}^{\alpha} \tilde{\lambda}^{\beta} \right]
  &= I^{\alpha\beta} + O_P (N^{-1})
  ,
  \label{formula:AEP_lambda2} 
\end{align}
where $\alpha, \beta, \gamma, \delta$ run through the indices $\{ 1, \cdots, p, \bar{1}, \cdots, \bar{p} \}$, and
\begin{align}
  L_{\alpha_1 \cdots \alpha_a} &:= \frac{1}{N} \pd{\alpha_1} \cdots \pd{\alpha_a} l_{\hat{\theta}}^{(N)}
  \label{formula:AE_L}, \\
  I^{\alpha_1 \cdots \alpha_a} &:= \int_{\sqrt{N}(\Theta - \hat{\theta})} \tilde{\lambda}^{\alpha_1} \cdots \tilde{\lambda}^{\alpha_a}
  \frac{1}{G^{(N)}} e^{- \frac{1}{2} J^{(N)}_{\alpha\beta} \tilde{\lambda}^{\alpha} \tilde{\lambda}^{\beta}} d\tilde{\lambda}
  \label{formula:AE_I}
\end{align}
for $\alpha_1, \cdots, \alpha_a \in \{1, \cdots, p, \bar{1}, \cdots, \bar{p} \}$.
Note that $L_{\alpha_1 \cdots \alpha_a}$ and $I^{\alpha_1 \cdots \alpha_a}$ are complex-valued random variables
because they depend on the realization of the observation $z^{(N)} \in \mathbb{C}^N$ from the process.

The expansion of (\ref{formula:AEP_S1}) now becomes
\begin{align}
  \hat{S}_{\pi}^{(N)} (\omega)
  &= S \,\bigl( \omega \mid \hat{\theta} \,\bigr)
  + \frac{1}{N} B^{(N)}_{\pi} \,\bigl( \omega \mid \hat{\theta} \,\bigr)
  + O_P \bigl( N^{-\frac{3}{2}} \bigr)
  \label{formula:AEP_S2}
  ,
\end{align}
where $B^{(N)}_{\pi}$ is the $O_P (1)$ term defined as
\begin{align}
  B^{(N)}_{\pi} \left( \omega \mid \theta \right) &:= \frac{1}{6} L_{\beta\gamma\delta} I^{\alpha\beta\gamma\delta} \left( \pd{\alpha} S \left( \omega \mid \theta \right) \right) 
    + I^{\alpha\beta} \left( \pd{\beta} \log \pi \bigl( \theta \bigr) \right) \left( \pd{\alpha} S \left( \omega \mid \theta \right) \right) 
    + \frac{1}{2} I^{\alpha\beta} \left( \pd{\alpha} \pd{\beta} S \left( \omega \mid \theta \right) \right)
    . 
\end{align}

Utilizing the complex-valued tensorial Hermite polynomials defined in Appendix \ref{section:hermite_polynomial}
and by Proposition \ref{prop:expectation_loglikehood} in Appendix \ref{proof:expectation_loglikehood},
we obtain an asymptotic expansion
$ B^{(N)}_{\pi} \left( \omega \mid \theta \right) = G^{(N)}_{\pi} \left( \omega \mid \theta \right) + H^{(N)} \left( \omega \mid \theta \right) + O_P(N^{-\frac{1}{2}}) $,
which yields the asymptotic expansion in (\ref{formula:AEP_S}).

Functions $G^{(N)}_{\pi}$ and $H^{(N)}$ represent the parallel and orthogonal parts of the quantity $N \bigl( \hat{S}^{(N)}_{\pi} - S_{\hat{\theta}^{(N)}} \bigr)$, respectively; see also \cite{komaki2006shrinkage, tanaka2011asymptotic}.

\section{Existence of Bayesian predictive power spectral densities for $AR(p; \mathbb{C})$}
\label{section:existenece_BPSD}

First, we prove that the $\kappa$-prior (\ref{definition:alpha_prior}) for $AR(p; \mathbb{C})$ is proper on $\Tilde{\Theta}_1 = U \times \cdots \times U$ if $\kappa < 1$ and is improper if $\kappa \geq 1$.
If $\kappa \leq 0$, then the $\kappa$-prior $\pi^{(\kappa)}$ is certainly integrable on $\Tilde{\Theta}_1$; therefore, we may assume $\kappa > 0$.
Because
\begin{align*}
| 1 - \xi^i \bar{\xi}^j |^2 - | \xi^i - \xi^j |^2 = ( 1 - |\xi^i|^2 )( 1 - |\xi^j|^2) \geq 0,
\end{align*}
we have $| \xi^i - \xi^j |^2 \,/\, | 1 - \xi^i \bar{\xi}^j |^2 \leq 1$.
Thus, for $0 < \kappa < 1$,
\begin{align*}
  \int_{\Tilde{\Theta}_1} \pi^{(\kappa)} (\xi) \,d\xi
  &= \int_{\Tilde{\Theta}_1} \frac{ \prod_{1 \leq i < j \leq q} | \xi^i - \xi^j |^2 }{ \prod_{i=1}^p \prod_{j=1}^p \bigl( 1 - \xi^i \bar{\xi}^j \bigr)^{\kappa} } \,d\xi \\
  &= \int_{\Tilde{\Theta}_1} \frac{1}{\prod_{i=1}^p \bigl( 1 - |\xi^i|^2 \bigr)^{\kappa}} \left( \prod_{1 \leq i < j \leq p} \frac{|\xi^i - \xi^j|^2}{|1 - \xi^i \bar{\xi}^j|^2} \right)^{\kappa} 
  \prod_{1 \leq i < j \leq p} |\xi^i - \xi^j|^{2(1 - \kappa)} \,d\xi \\
  &\leq 2^{p(p-1)(1-\kappa)} \prod_{i=1}^p \left( \int_U \frac{1}{\bigl( 1 - |\xi^i|^2 \bigr)^{\kappa}} d\xi^i \right) \\
  &= 2^{p(p-1)(1-\kappa)} \left( \frac{\pi}{1 - \kappa} \right)^p
\end{align*}
because $\int_U  \bigl( 1 - |\xi^i|^2 \bigr)^{-\kappa} d\xi^i = \int_{0}^1 \int_{-\pi}^{\pi} \bigl( 1 - r^2 \bigr)^{-\kappa} \,r \,dr \,d\theta = \frac{\pi}{1 - \kappa}$.
Therefore, the $\kappa$-prior $\pi^{(\kappa)}$ is integrable on $\Tilde{\Theta}_1$ if $\kappa < 1$.

Set
\begin{align*}
  m := \min_{\, \xi \,\in\, \Xi} \prod_{1 \leq i < j \leq p} \frac{|\xi^i - \xi^j|^2}{|1 - \xi^i \bar{\xi}^j|^{2\kappa}} > 0,
\end{align*}
where $\Xi := V_1 \times \cdots \times V_p \subset \Tilde{\Theta}_1$ and
\begin{align*}
  V_i := \left\{ \xi \in U \;\middle|\; \frac{1}{2} < |\xi| < 1, \frac{2\pi}{N} (i - 1/2) < \arg \xi < \frac{2\pi}{N} i \right\}.
\end{align*}
We see that $\pi^{(\kappa)}$ is not integrable on $\Tilde{\Theta}_1$ if $\kappa \geq 1$, because
\begin{align*}
  \int_{\Tilde{\Theta}_1} \pi^{(\kappa)} (\xi) \,d\xi &\geq \int_{\Xi} \pi^{(\kappa)} (\xi) \,d\xi 
  \geq m \prod_{i=1}^p \int_{V_i} (1 - |\xi|^2)^{-\kappa} \,d\xi = +\infty.
\end{align*}

Next, we prove that a function $p^{(N)}_{\xi} \bigl( z^{(N)} \bigr) \, \pi^{(\kappa)} \bigl( \xi \bigr)$ of $\xi$ is integrable on the parameter space $\Tilde{\Theta}_1$ if $\kappa < 2$.
The explicit form of the determinant $\bigl| \Sigma^{(N)} \bigr|$ of the variance-covariance matrix $\Sigma^{(N)}$ of $AR(p; \mathbb{C})$ of the form in (\ref{definition:parameter_xi}) is
$\gamma_0 = \bigl| \Sigma^{(1)} \bigr| \leq \cdots \leq \bigl| \Sigma^{(p)} \bigr| = \bigl| \Sigma^{(p + 1)} \bigr| = \cdots = \prod_{i = 1}^p \prod_{j = 1}^p \bigl( 1 - \xi^i \bar{\xi}^j \bigr)^{-1}$;
see Section 5.5 (b), (c), and (d) in \cite{grenander1958toeplitz} or Theorem 3.1 in \cite{shaman1976approximations}.
Thus, if $N \geq p$, we have
\begin{align*}
  p^{(N)} \bigl( z^{(N)} \mid \xi \bigr) \, \pi^{(\kappa)} \bigl( \xi \bigr)
  &=
  \pi^{-N} \,
  \bigl| \Sigma^{(N)} \bigr|^{-1}
  \,
  e^{ - z^{(N)*} \bigl( \Sigma^{(N)} \bigr)^{-1} z^{(N)} }
  \frac{ \prod_{1 \leq i < j \leq q} | \xi^i - \xi^j |^2 }{ \prod_{i=1}^p \prod_{j=1}^p \bigl( 1 - \xi^i \bar{\xi}^j \bigr)^{\kappa} }
  \\
  &\leq \pi^{-N}  \frac{ \prod_{1 \leq i < j \leq q} | \xi^i - \xi^j |^2 }{ \prod_{i=1}^p \prod_{j=1}^p \bigl( 1 - \xi^i \bar{\xi}^j \bigr)^{\kappa - 1} } \\
  &= \pi^{-N} \bigl( \pi^{(\kappa - 1)} (\xi) \bigr)
\end{align*}
for $z^{(N)} \in \mathbb{C}^N$.
If $\kappa < 2$ and $N \geq p$, then $\int_{\Tilde{\Theta}_1} p^{(N)}_{\xi} \bigl( z^{(N)} \bigr) \, \pi^{(\kappa)} \bigl( \xi \bigr) \,d\xi$ is bounded, regardless of the sample $z^{(N)} \in \mathbb{C}^N$.
Therefore, the Bayesian predictive power spectral densities $\hat{S}^{(N)}_{\pi^{(\kappa)}}$ for $AR(p; \mathbb{C})$ based on the $\kappa$-prior $\pi^{(\kappa)}$ exists if $\kappa < 2$ and $N \geq p$.

The posterior $\pi^{(\kappa)} \bigl( \xi \mid z^{(N)} \bigr)$ given an observation $z^{(N)}$ based on the prior $\pi^{(\kappa)} ( \xi )$ is calculated as
\begin{align}
  \pi^{(\kappa)} \bigl( \xi \mid z^{(N)} \bigr)
  \propto p^{(N)} \bigl( z^{(N)} \mid \xi \bigr) \, \pi^{(\kappa)} \bigl( \xi \bigr)
  \propto \pi^{(\kappa - 1)} (\xi) \times e^{ - z^{(N)*} \bigl( \Sigma^{(N)} \bigr)^{-1} z^{(N)} }
\end{align}
if $N \geq p$, which shows that the family of $\kappa$-priors is not closed under sampling, because of the exponential term.
The conjugation property of the prior family is a topic for future research.

\section{Proof of the Main Theorem}
\label{proof:main}

Here, we prove the Main Theorem.
The generalization, including the i.i.d.\ case, for this theorem is discussed in Section \ref{section:generalization}. 
Except Proposition \ref{proposition:AEP_R_DIFF_ALPHA}, most of the propositions presented in this section are merely complexified versions of previous works; see \cite{komaki2006shrinkage, tanaka2011asymptotic}.
In particular, the proofs for Proposition \ref{proposition:AEP_R} and Proposition \ref{proposition:AEP_R_DIFF} only rely on the form (\ref{formula:AEP_S}) of the asymptotic expansion of Bayesian predictive power spectral densities (or Bayesian predictive distributions for the i.i.d.\ case). 
Therefore, we present the minimal outlines for these proofs because the procedures in previous works can be directly applied to them.

We first present the asymptotic expansion of the risk $R \bigl( \hat{S}^{(N)}_{\pi} \mid \theta \bigr)$ for a Bayesian predictive power spectral density $\hat{S}^{(N)}_{\pi}$.
The original proof of Proposition \ref{proposition:AEP_R} for the i.i.d.\ case was introduced in \cite{komaki2006shrinkage},
and the proof for real-valued processes is reported in \cite{tanaka2011asymptotic}.

The covariant derivative $\overset{\tiny{(e)}\;\;}{\nabla_\alpha}$ of the vector field $V^{\beta}$ is defined as
$ \overset{\tiny{(e)}\;\;}{\nabla_\alpha} V^{\beta} := \pd{\alpha} V^{\beta} + \Guddu{(e)\;\;}{\alpha}{\gamma}{\beta} V^{\gamma} $,
where $\Guddu{(e)\;\;}{\alpha}{\gamma}{\beta} := \Guddu{(m)\;\;}{\alpha}{\gamma}{\beta} - \Tddd{\alpha}{\gamma}{\delta} \, \guu{\delta}{\beta}$.

\begin{prop}
  \label{proposition:AEP_R}
  For a complex parameter space $\Theta$ and a possibly improper prior $\pi$, we have
  \begin{align}
    R \bigl( \hat{S}^{(N)}_{\pi} \mid \theta \bigr)
    = \frac{1}{2 N^2} \guu{\alpha}{\beta} \left( \pd{\alpha} \log \frac{\pi}{\pi_J} + \frac{1}{2} \Td{\alpha} \right) \left( \pd{\beta} \log \frac{\pi}{\pi_J} + \frac{1}{2} \Td{\beta} \right)
    + \frac{1}{N^2} \overset{\tiny{(e)}\;\;}{\nabla_\alpha} \left( \guu{\alpha}{\beta} \left(  \pd{\beta} \log \frac{\pi}{\pi_J} + \Td{\beta} \right) \right)
    + C + O(N^{-\frac{5}{2}}) \label{formula:AEP_R}, 
  \end{align}
  where $C$ is a term independent of the prior $\pi$.
\end{prop}

\begin{IEEEproof}
  Recall that the asymptotic expansion of a Bayesian predictive power spectral density $\hat{S}^{(N)}_{\pi}$ is given by (\ref{formula:AEP_S}).
  We follow the procedure in \cite{tanaka2011asymptotic} but replace the summation rule $i, j, k = 1, \cdots, p$ with $\alpha, \beta, \gamma = 1, \cdots, p, \bar{1}, \cdots, \bar{p}$.
\end{IEEEproof}

The asymptotic expansion of the risk difference $R \bigl( \hat{S}^{(N)}_{\pi_1} \mid \theta \bigr) - R \bigl( \hat{S}^{(N)}_{\pi_2} \mid \theta \bigr)$ is calculated as follows;
see also \cite{komaki2006shrinkage, tanaka2011asymptotic}.

\begin{prop}
  \label{proposition:AEP_R_DIFF}
  For a complex parameter space $\Theta$ and a possibly improper prior $\pi$, we have
  \begin{align*}
    N^2 \left( R \bigl( \hat{S}^{(N)}_{\pi_J} \mid \theta \bigr) - R \bigl( \hat{S}^{(N)}_{\pi} \mid \theta \bigr) \right) 
    &= \guu{\alpha}{\beta}
      \left( \pd{\alpha} \log \frac{\pi}{\pi_J} \right) \left( \pd{\beta} \log \frac{\pi}{\pi_J} \right)
    - \left( \frac{\pi}{\pi_J} \right)^{-1} \Delta \left( \frac{\pi}{\pi_J} \right) 
    + O (N^{-\frac{1}{2}}) \\
    &= 
    -2 \left( \frac{\pi}{\pi_J} \right)^{-\frac{1}{2}} \Delta \left( \frac{\pi}{\pi_J} \right)^{\frac{1}{2}}
    + O (N^{-\frac{1}{2}}) \nonumber
  \end{align*}
  for $\theta \in \Theta$.
\end{prop}

\begin{IEEEproof}
  Use (\ref{formula:AEP_R}) and follow \cite{tanaka2011asymptotic}.
\end{IEEEproof}

\begin{prop}
  \label{proposition:AEP_R_DIFF_ALPHA}
  Let $\phi$ be a positive continuous function on a K\"ahler parameter space $\Theta$, and define a family $\pi^{(\kappa)} := \phi^{-\kappa+1} \pi_J$ of priors for $\kappa \in \mathbb{R}$.
  Then, we have
  \begin{align*}
    N^2 \left( R \bigl( \hat{S}^{(N)}_{\pi_1} \mid \theta \bigr) - R \bigl( \hat{S}^{(N)}_{\pi_2} \mid \theta \bigr) \right) 
    = - (\kappa_1 - \kappa_2) \frac{\Delta \phi}{\phi} 
    + (\kappa_1 - \kappa_2) (\kappa_1 + \kappa_2) \gHuAu{i}{j} \bigl( \pHd{i} \log \phi \bigr) \bigl( \pAd{j} \log \phi \bigr) + O (N^{-\frac{1}{2}})
  \end{align*}
  for $\theta \in \Theta$, where $\pi_1 := \pi^{(\kappa_1)}$ and $\pi_2 := \pi^{(\kappa_2)}$.
\end{prop}

\begin{IEEEproof}
  Let $\kappa_1 = 1$ and $\kappa_2 = \kappa$.
  Using (\ref{formula:laplacian_kappa}) for $a = \frac{1-\kappa}{2}$, we have
  \begin{align*}
    N^2 \left( R \bigl( \hat{S}^{(N)}_{\pi_J} \mid \theta \bigr) - R \bigl( \hat{S}^{(N)}_{\pi^{(\kappa)}} \mid \theta \bigr) \right)
    &= - 2 \, \frac{ \Delta \phi^{ \frac{1-\kappa}{2}} }{ \phi^{ \frac{1-\kappa}{2}}} + O (N^{-\frac{1}{2}}) \\
    &= (1 - \kappa) \frac{\Delta \phi}{\phi} + (\kappa^2 - 1) \gHuAu{i}{j} \bigl( \pHd{i} \log \phi \bigr) \bigl( \pAd{j} \log \phi \bigr)
    + O (N^{-\frac{1}{2}})
    .
  \end{align*}
  The formula
  \begin{align*}
    N^2 \left( R \bigl( \hat{S}^{(N)}_{\pi_1} \mid \theta \bigr) - R \bigl( \hat{S}^{(N)}_{\pi_2} \mid \theta \bigr) \right)
    = - N^2 \left( R \bigl( \hat{S}^{(N)}_{\pi_J} \mid \theta \bigr) - R \bigl( \hat{S}^{(N)}_{\pi_1} \mid \theta \bigr) \right) 
    - N^2 \left( R \bigl( \hat{S}^{(N)}_{\pi_J} \mid \theta \bigr) - R \bigl( \hat{S}^{(N)}_{\pi_2} \mid \theta \bigr) \right)
  \end{align*}
  yields the statement.
\end{IEEEproof}

If there exists a positive continuous eigenfunction $\phi > 0$ of the Laplacian $\Delta$ with a negative eigenvalue $-K < 0$, then $\frac{\Delta \phi}{\phi} = -K$,
which yields the proof of Theorem \ref{theorem:main}.

\section{Relation with $\alpha$-parallel priors}
\label{section:alpha_parallel}

In this appendix, we show the relation between $\kappa$-priors and $\alpha$-parallel priors.

First, we define the $\alpha$-connection.
Let $\Theta$ be a complex parameter space in $\mathbb{C}^p$.
For $\alpha \in \mathbb{R}$, set
\begin{align}
  \Guddd{(\alpha)}{\beta}{\gamma}{\delta} &:= \Guddd{(m)}{\beta}{\gamma}{\delta} - \frac{1 + \alpha}{2} \, \Tddd{\beta}{\gamma}{\delta},
  \label{definition:alpha_connection} \\
  \Guddu{(\alpha)}{\beta}{\gamma}{\delta} &:= \Guddd{(\alpha)}{\beta}{\gamma}{\epsilon} \, \guu{\epsilon}{\delta}
\end{align}
for $\beta, \gamma, \delta \in \{ 1, \cdots, p, \bar{1}, \cdots, \bar{p} \}$.
For complex-valued stationary ARMA processes, the quantities $\gdd{\alpha}{\beta}$, $\Tddd{\alpha}{\beta}{\gamma}$, and $\Guddd{(m)}{\alpha}{\beta}{\gamma}$ on the right-hand side of the equations are defined as (\ref{definition:M_g}), (\ref{definition:M_T}), and (\ref{definition:M_Gm}), respectively.
For the i.i.d.\ case, these quantities are defined as (\ref{definition:K_g}), (\ref{definition:K_T}), and (\ref{definition:K_Gm}), respectively.
Note that $\Guddd{(-1)}{\beta}{\gamma}{\delta}$ corresponds to $\Guddd{(m)}{\beta}{\gamma}{\delta}$.
On the other hand, $\Guddd{(+1)}{\beta}{\gamma}{\delta} = \Guddd{(m)}{\beta}{\gamma}{\delta} - \Tddd{\beta}{\gamma}{\delta}$ is often denoted by $\Guddd{(e)}{\beta}{\gamma}{\delta}$.
The quantity $\Guddu{(\alpha)}{\beta}{\gamma}{\delta}$ is called the coefficients of the $\alpha$-connection $\nabla^{(\alpha)}$ on the parameter space.
The $\alpha$-connection with $\alpha = -1$ is called the mixture connection ($m$-connection) $\nabla^{(m)}$, and the $\alpha$-connection with $\alpha = +1$ is called the exponential connection ($e$-connection) $\nabla^{(e)}$.
The $\alpha$-connection with $\alpha = 0$ is called the Levi--Civita connection $\nabla^{(0)}$.
From (\ref{definition:alpha_connection}), we have
\begin{align}
  \Guddd{(\alpha)}{\beta}{\gamma}{\delta} = \Guddd{(0)}{\beta}{\gamma}{\delta} - \frac{\alpha}{2} \, \Tddd{\beta}{\gamma}{\delta}
  \,.
\end{align}
For the geometrical interpretation of the $\alpha$-connection $\nabla^{(\alpha)}$ on the statistical manifold, see \cite{amari1983differential, amari1985differential}.

If the complex parameter space $\Theta$ is K\"ahler, the nontrivial elements of the coefficients $\Guddu{(0)\;}{\beta}{\gamma}{\delta}$ of the $0$-connection $\nabla^{(0)}$ are only $\GuHdHdHu{(0)\;}{i}{j}{k}$ and $\GuAdAdAu{(0)\;}{i}{j}{k}$;
other coefficients, such as $\Guddu{(0)\;}{i}{j}{\bar{k}}$, $\Guddu{(0)\;}{i}{\bar{j}}{k}$, and $\Guddu{(0)\;}{i}{\bar{j}}{\bar{k}}$, are all equal to $0$; see Section 8.5 in \cite{nakahara2003geometry}.

The non-negative function $\rho^{(\alpha)}$ with $\alpha \in \mathbb{R}$ defined on the complex parameter space $\Theta$ is called an $\alpha$-parallel prior, if
\begin{align} 
  \pd{\beta} \, \rho^{(\alpha)} = \Guddu{(\alpha)}{\beta}{\gamma}{\gamma} \, \rho^{(\alpha)}
\end{align}
for $\beta, \gamma \in \{ 1, \cdots, p, \bar{1}, \cdots, \bar{p} \}$ on the complex parameter space $\Theta$; see \cite{takeuchi2005spl}.
The quantity $\left( \pd{\beta} - \Guddu{(\alpha)}{\beta}{\gamma}{\gamma} \right) \rho^{(\alpha)}$ is referred to as a contravariant derivative of $\rho^{(\alpha)}$ with respect to the $\alpha$-connection $\nabla^{(\alpha)}$.
The $\alpha$-parallel prior $\rho^{(\alpha)}$ corresponds to the volume element that is parallel with respect to the $\alpha$-connection $\nabla^{(\alpha)}$.
In general, the existence of the $\alpha$-parallel prior for each $\alpha \in \mathbb{R}$ is not guaranteed.
However, we can provide a sufficient condition for the existence of a family of $\alpha$-parallel priors on a K\"ahler parameter space.

Let us compute the contravariant derivatives of priors with respect to the $\alpha$-connection $\nabla^{(\alpha)}$.
Let $\phi$ be a non-negative function globally defined on a K\"ahler parameter space $\Theta \subset \mathbb{C}^p$.
We define a family $\{ \pi^{(\kappa)} \}_{\kappa \in \mathbb{R}}$ of priors by $\pi^{(\kappa)} := \phi^{-\kappa+1} \pi_J$,
where $\pi_J$ denotes the Jeffreys prior.
Note that $\pi^{(1)}$ corresponds to the Jeffreys prior $\pi_J$.
As $\Theta$ is K\"ahler, we have
\begin{align}
  \GuHd{(0)}{i} &:= \GuHdHdHu{(0)}{i}{j}{j} + \GuHdAdAu{(0)}{i}{k}{k} = \GuHdHdHu{(0)}{i}{j}{j} = \GuHdHdAd{(0)}{i}{j}{k} \gHuAu{j}{k}
  \nonumber\\
  &= \bigl( \pHd{i} \gHdAd{j}{k} \bigr) \gHuAu{j}{k} = \pHd{i} \log \pi_J,
\end{align}
where we used the Jacobi formula (\ref{formula:jacobi}).
As
\begin{align}
  \GuHd{(0)}{i} \pi^{(\kappa)} = \bigl( \pHd{i} \log \pi_J \bigr) \phi^{-\kappa+1} \pi_J = \phi^{-\kappa+1} \, \pHd{i} \pi_J,
\end{align}
we have
\begin{align}
  \GuHd{(\alpha)}{i} \pi^{(\kappa)} &= \left( \GuHd{(0)}{i} - \frac{\alpha}{2} \THd{i} \right) \pi^{(\kappa)}
  \nonumber \\
  &= \phi^{-\kappa+1} \, \pHd{i} \pi_J - \frac{\alpha}{2} \, \THd{i} \, \phi^{-\kappa+1} \pi_J,
\end{align}
where $\THd{i} := \THdHdAd{i}{k}{j} \gHuAu{k}{j} + \THdAdHd{i}{j}{k} \gHuAu{k}{j}$.
On the other hand, we have
\begin{align}
  \pHd{i} \pi^{(\kappa)} &= \phi^{-\kappa+1} \, \pHd{i} \pi_J - (\kappa - 1) \phi^{-\kappa+1} (\pHd{i} \log \phi) \pi_J .
\end{align}
Thus, the contravariant derivative of $\pi^{(\kappa)}$ with respect to the $\alpha$-connection $\nabla^{(\alpha)}$ is
\begin{align}
  \left( \pHd{i} - \GuHd{(\alpha)}{i} \right) \pi^{(\kappa)}
  = \left( \frac{\alpha}{2} \THd{i} - (\kappa - 1) (\pHd{i} \log \phi) \right) \phi^{-\kappa+1} \pi_J \,.
\end{align}
We see that if $\pHd{i} \log \phi$ is proportional to $\THd{i}$, the proposed prior $\pi^{(\kappa)}$ is an $\alpha$-parallel prior for some $\alpha$.
Thus, we have the following proposition.
\begin{prop}
  If $\THd{i} = -4 \, c \, \pHd{i} \log \phi$ for some $c \neq 0$, then
  the prior $\pi^{(\kappa)}$ is an $\alpha$-parallel prior with $\alpha = (1 - \kappa) / 2c$.
\end{prop}

Moreover, the term $\left( \pd{\alpha} \log \frac{\pi}{\pi_J}+ \frac{1}{2} T_{\alpha} \right)$, which appears in the definition (\ref{definiton:G}) of the parallel part $G^{(N)}_{\pi}$ of the risk difference between the Bayesian predictive distribution and the estimative distribution with the maximum likelihood estimator, now becomes
\begin{align}
  \pHd{i} \, \log \frac{\pi^{(\kappa)}}{\pi_J}+ \frac{1}{2} \, \THd{i} = - \bigl( \kappa - (1 - 2c) \bigr) \, \pHd{i} \log \phi
\end{align}
for $i = 1, \cdots, p$.
Thus, we have the following proposition.
\begin{prop}
  If $\THd{i} = -4 \, c \, \pHd{i} \log \phi$ for some $c \neq 0$, then
  the risk of the Bayesian predictive distribution asymptotically dominates the risk of the estimative distribution with the maximum likelihood estimator when $\kappa = 1 - 2 c$, i.e., $\alpha = (1 - \kappa) / 2c$ = 1.
\end{prop}

For the non-negative function $\phi$ defined as (\ref{definition:phi}) for $\mathrm{AR} (p; \mathbb{C})$,
we have $\THd{i} = -4 \pHd{i} \log \phi = 4 \bar{\xi}^j \gHdAd{i}{j}$, i.e., $c = 1$.
Thus, $\pi^{(\kappa)}$ is an $\alpha$-parallel prior with $\alpha = (1 - \kappa) / 2$ for $\mathrm{AR} (p; \mathbb{C})$.
In particular, the proposed prior $\pi^{(-1)}$, which asymptotically achieves the constant risk improvement, is a ($+1$)-parallel prior.
Moreover, the risk of the Bayesian predictive power spectral density based on the proposed prior $\pi^{(-1)}$ asymptotically dominates the risk of the estimative predictive power spectral density with the maximum likelihood estimator.

\section*{Acknowledgment}

We express our gratitude to Keisuke Yano for helpful comments.

\ifCLASSOPTIONcaptionsoff
  \newpage
\fi

\vfill

\end{document}